\documentclass[12pt]{article}
\usepackage{amsfonts,amssymb,version,amsmath}
\def\A{{\mathcal A}}
\def\C{\mathbb C}
\def\D{{\mathcal D}}
\def\EE{{\mathbb E}}
\def\F{{\mathcal F}}
\def\I{{\mathcal I}}

\def\M{{\mathcal M}}
\def\OO{{\mathcal O}}

\def\R{\mathbb R}

\def\Z{{\mathbb Z}}

\def\deg{\mathop{\rm deg}\nolimits} 
\def\diag{\mathop{\rm diag}\nolimits}

\def\id{\mathop{\rm id}\nolimits} 
 
\def\ker{\mathop{\rm ker}\nolimits}

\mathchardef\mhyphen="2D

\def\rank{\mathop{\rm rank}\nolimits}
\def\res{\mathop{\rm res}\nolimits}
\def\trace{\mathop{\rm trace}\nolimits}

\let\da\downarrow

\let\hra\hookrightarrow

\let\ov\overline

\let\un\underline

\let\wh\widehat

\let\pa\partial

\def\stm{\refstepcounter{theorem}\paragraph{\thetheorem}}

\def\proof{\paragraph{Proof}}

\makeatletter \def\l@section{\@dottedtocline{1}{0em}{1.2em}} \makeatother

\begin{document}

\centerline{\Large\bf The Narasimhan-Seshadri Theorem revisited}

\bigskip

\centerline{\bf Nitin Nitsure}

\bigskip

\centerline{\small MSC2020: 14F10, 14H60}

\begin{abstract} Let $X$ be a compact Riemann surface. 
  The famous Narasimhan-Seshadri theorem [13] of 1965 uses the
  Grothendieck construction [4] of 1956 that associates
  vector bundles $\EE(\sigma)$ on $X$ to representations $\sigma$
  of a certain Fuchsian group $\pi$. Narasimhan and Seshadri
  show that by taking the representations $\sigma$ to be irreducible unitary of a certain kind, this
  exactly gives all stable vector bundles on $X$ of a given rank and degree.

  In this note we reformulate the correspondence from representations to bundles, which leads to
  simpler statements and proofs. The Fuchsian group $\pi$ is replaced by the punctured fundamental group 
  $\pi_1(X-x)$ where $x\in X$. The Grothendieck bundles $\EE(\sigma)$ then become 
  Deligne's logarithmic extensions to $X$ of bundles with connections on $X-x$ associated to
  representations of $\pi_1(X-x)$ with scalar local monodromy. We also report how some 
  ideas from algebraic geometry (which were all in place by 1970) have 
  simplified some aspects of the original proof over the decades. 

  This simplified approach works equally well for all values of the genus $g$, removing the
  restriction $g\ge 2$ in the 1965 original. 
  Finally, we comment that such a logarithmic reformulation extends to related kinds of bundles
  such as principal bundles with reductive structure groups or parabolic bundles. 
\end{abstract}

{\footnotesize

\tableofcontents

}

\newpage


\pagestyle{myheadings}
\markright{Narasimhan-Seshadri Theorem revisited.}


\section{Introduction}
The statement and the proof of the Narasimhan-Seshadri theorem in their famous 1965 paper [13]
can be seen to consist of the following five main steps.
We begin by describing these five steps in broad terms, and then say how they can be modified and
simplified. The principal modification is to replace Grothendieck's 1956 construction of bundles
by Deligne's 1970 logarithmic extensions of connections. We also mention some other modifications
which have become commonplace over time, involving ideas such as Quot schemes,
algebraic deformation theory,
Geometric Invariant Theory, good quotients, projective moduli variety of semistable bundles.
These ideas were all in place by 1970, but were not available to Narasimhan and Seshadri in 1965.

{\bf (1) The group $\pi$.} (Poincar\'e.) Let $X$ be a compact Riemann surface of genus $g\ge 1$.
For any integer $N\ge 1$, Poincar\'e showed that there is
a simply connected Riemann surface $Y$ and a ramified Galois cover $p: Y\to X$, which is
ramified at most over a single point $x\in X$, such that the inertia subgroup of $\pi = Aut(Y/X)$
at any point of the fiber $p^{-1}(x)$ is cyclic of order $N$. If $N=1$ then $p$ is an unramified covering,
and $\pi$ becomes $\pi_1(X)$. 

{\bf (2) Grothendieck construction.} Let $\sigma: \pi \to GL(n)$ be any representation
with a certain prescribed behaviour on the inertia subgroups of $\pi$. 
Grothendieck considers in [4] 
the trivial vector bundle $\OO_Y^n$ on $Y$ together with a $\pi$-action via $\sigma$, lifting
the $\pi$-action on $Y$. He then takes its 
invariant direct image sheaf $\EE(\sigma)$ on $X$, which is a locally free coherent $\OO_X$-module.
Grothendieck showed that if we took $N$ to be a multiple of $n$, and chose suitably 
the restriction of $\sigma$ to any inertia subgroup then the bundle $\EE(\sigma)$ will have
a prescribed degree $0\le m < n$. 
The Grothendieck construction of  $\EE(\sigma)$ from $\sigma$ 
was somehow implicit in the 1938 paper of Weil [18], which Grothendieck reformulated
in bundle theoretic and sheaf theoretic terms in 1956.

{\bf (3) Irreducible unitary representations and stability.}
Narasimhan and Seshadri were intrigued by Weil's
1938 paper, in particular, by Weil's wondering about the significance of unitary representations.
In [13], they apply the above Grothendieck construction to unitary representations of $\pi$. 
They showed that if $\sigma$ is an irreducible unitary representation with prescribed behaviour
on isotropies, then the bundle $\EE(\sigma)$ is stable, and moreover, this gives an injective map
from the set of all conjugacy classes of representations to the set of isomorphism classes of
stable bundles. The proof of stability of $\EE(\sigma)$ and injectivity is ultimately
based on the classical theory of plurisubharmonic functions on open disks. 

{\bf (4) Moduli space of stable bundles.} We must remember that the work on [13] was done 
based on the mathematics of 1950's (or older), but without any input from the ongoing revolution
in Algebraic Geometry at the hands of Grothendieck and his school. The only modern material
they used was the definition of a stable bundle by Mumford, without any further 
input from the ongoing development of GIT. 
Instead, the paper [13]
uses the analytic deformation theory of Kodaira and Spencer [5] to construct a holomorphic
manifold $S_s(n,m)$ of stable vector bundles of rank $n$ degree $m$ on $X$. They also
show that $S_s(n,m)$ is connected.


{\bf (5) Open embedding and surjectivity.} To prove that every stable bundle $E$ is of the
form $\EE(\sigma)$, Narasimhan and Seshadri use what they call the `Poincar\'e continuity method'. 
They show that the map from the space of conjugacy classes of representations to the moduli of 
stable bundles is an open embedding and a proper map, and so it is an analytic isomorphism
as $S_s(n,m)$ is connected.

In our reformulation of the above material, the overall steps are similar, but there are changes
inside each step. These are broadly as follows.

{\bf (1*) The group $\pi_1(X-x,y)$.} In our reformulation, we have no use for the
ramified cover of Poincar\'e. Instead, we just replace $\pi$ in (1) by the fundamental
group $\pi_1(X-x,y)$ of the punctured Riemann surface, which is a free group
on $2g$ generators where $g$ is the genus of $X$. Moreover, we do not have any need
for the integer $N$ of (1), as this new $\pi$ will give us bundles of all ranks and degrees.
The further reformulated steps uniformly work for all values $g\ge 0$ of the genus.

{\bf (2*) The Deligne construction.} We consider representations $\rho : \pi_1(X-x,y)\to GL(n)$
for which the loop $c\in \pi_1(X-x,y)$ around the puncture $x$ is mapped to
$e^{2\pi i m/n} I$ where $m,n \in \Z$ with $n\ge 1$.
To such a representation there is naturally associated a holomorphic vector bundle $E_{\rho}$ on $X-x$
with a holomorphic connection $\nabla_{\rho}$. The Deligne construction of logarithmic extensions [3] 
attaches to this data an extension $E(\rho)$ of $E_{\rho}$ to $X$, which is a holomorphic vector bundle
of rank $n$, degree $m$ on $X$, such that 
$\nabla_{\rho}$ extends to a logarithmic connection on $E(\rho)$, which has residue $(-m/n)I$ at $x$.
The bundle $E(\rho)$ is naturally isomorphic to the Grothendieck bundle $\EE(\sigma)$
where $\sigma : \pi \to GL(n)$ is induced by $\rho$. 
This replaces the Grothendieck construction in our reformulation.

{\bf (3*) Irreducible unitary representations and stability.} By arguments very similar
to those in [13], again based on the behaviour of plurisubharmonic functions, we show that 
if $\rho$ is an irreducible unitary representation with $\rho(c) = e^{2\pi i m/n} I$,
then the bundle $E(\rho)$ is stable, and moreover, this gives an injective map
from the set of all conjugacy classes of such representations to the set of isomorphism classes of
stable bundles.

{\bf (4*) Moduli space of stable bundles.} We simply appeal to the modern 
GIT-based construction of the moduli of stable bundles, which allows us to completely
bypass the Kodaira-Spencer theory.

{\bf (5*) Open embedding and surjectivity.} This is almost the same as (5) except that
the local deformation
theory can now be done much better in the style of Grothendieck-Schlessinger-Artin. Another
simplification is to use the Mumford-Seshadri good moduli space $M^{ss}(X,n,m)$ of semistable
bundles to compactify the moduli $M^s(X,n,m)$, which makes the properness of the correspondence
obvious.

From about 1970, all treatments of the Narasimhan-Seshadri theorem have used algebraic deformation
theory and GIT-based moduli of stable and semistable bundles, so (4*) and (5*) are routine
for the past five decades or more. Another simplification due to Ramanan (which he discovered in
the 1960s, soon after the theorem was published)
is to use an argument based on decomposable multivectors in (3), and the same works for (3*).

This note is arranged as follows.


Section 2 sets up the basic definitions and conventions
that we will use through out. 
Section 3 rapidly recalls the well-known material on holomorphic ODEs, logarithmic connections,
and Deligne extensions in dimension $1$ that we need.
Section 4 treats the case of line bundles.
Section 5 uses the above machinery to state and prove our reformulation of Weil's theorem
on indecomposable bundles. 
Section 6 states the Main Theorem (Theorem \ref{main theorem}), which is
our reformulation of the Narasimhan-Seshadri Theorem.
Section 7 describes the precise relationship between the Grothendieck construction
used in [13] and the Deligne construction that we use.
Section 8 says how our formulation works for $g=0$ and $g=1$, the two cases that were
left out of [13].
Section 9 first recalls some well known facts about plurisubharmonic functions, and then
proves the statement (1) of the Main Theorem \ref{main theorem}.
Section 10 is devoted to a rapid sketch of the steps (4*) and (5*).

{\bf Historical note:} I had outlined this reformulation to 
Professors Narasimhan and Seshadri during the 2015 conference at CMI to mark the 50th anniversary
of the theorem. Both of them asked me to write it up as an article.
Here it is, ten years later. To my regret, I did not get around to it in their lifetimes.

\section{Notations and conventions}
Throughout, $X$ will denote a compact Riemann surface of arbitrary genus $g\ge 0$,
and $x,y\in X$ two distinct chosen points.
Let $r>1$, $\epsilon>0$ and let $(U_{r+\epsilon},z)$ be a holomorphic coordinate chart on
$X$ where $z(U_{r+\epsilon})$ is the open disk $|z|< r+\epsilon$ in $\C$
and such that $x$ and $y$ lie in $U_{r+\epsilon}$ and correspond to $z=0$ and $z=1$.
We denote by $U\subset U_{r+\epsilon}$ the open subset where $|z|< r$. 
We will denote by $c\in \pi_1(X-x,y)$ the element corresponding to the simple positive loop around $x$,
defined by the map $[0,1]\to U-x : t\mapsto e^{2\pi i t}$. 
As usual, a {\bf frame} over $y$ for
a local system or a vector bundle on $X$ will mean a basis $y^*$ for its fiber over $y$.

Following [3], we compose loops (or paths) on the left in defining fundamental groups
(or fundamental groupoids). For a pointed universal cover of a space with a base point,
the Galois group (deck transformation group) will act on the left on the total space,
and it is therefore the opposite of the fundamental group, which will act on the right. 

We assume a basic familiarity with the notions of local systems, holomorphic connections,
logarithmic connections, monodromy representations etc., as in Deligne [3]. An exposition of 
these ideas is available in [6]. We recall the following basic definitions for convenience.

\stm {\bf Logarithmic connections on Riemann surfaces and their residues.}
A logarithmic connection on $(X,x)$ is a pair $(E,\nabla)$, which
consists of a holomorphic vector bundle $E$ on $X$ equipped with a $\C$-linear homomorphism of
sheaves
$$\nabla : E \to \Omega_X^1(\log x)\otimes_{\OO_X}E$$
that satisfies the Leibniz rule $\nabla(f v) = df\otimes v + f \nabla(v)$ for
any local sections $f$ and $v$ of $\OO_X$ and $E$, respectively.
The line bundle $\Omega_X^1(\log x)$ of logarithmic differentials
%
%
on $(X,x)$ has local free
basis $dz/z$ over $U$. 
Note that the restriction $(E|_{X-x}, \nabla|_{X-x})$ is a holomorphic connection on $X-x$. 
The {\bf residue} of $(E,\nabla)$ at $x$ is the element $\res_x(\nabla) \in End(E_x)$
defined by the composite homomorphism
$$E \stackrel{\nabla}{\to} \Omega^1_X(\log x)\otimes_{\OO_X} E
\stackrel{P\otimes \id_E}{\to} \OO_x \otimes E = E_x$$
(which is $\OO_X$-linear, even though $\nabla$ is not so)
where $P: \Omega^1_X(\log x)\to \OO_x$ denotes the Poincar\'e 
residue map, with $P($dz/z$) = 1$. A logarithmic connection on $(X,x)$ is a holomorphic connection
on $X$ if and only if its residue is $0$.

We introduce various sets that we use later, by means of the following table.

\bigskip

{\small 
\begin{tabular}{ll}
{\bf $~~~~~~~$ Set} & {\bf $~~~~~~~~~~~~~~~$ Elements} \\


$Rep(X-x,y, n, \tau)$ & representations $\rho : \pi_1(X-x,y) \to GL(n)$ such that\\
                      & $\rho(c) = \tau I$ where $\tau \in \C^{\times}$. \\


$Rep(X-x, n, \tau)$ & conjugacy classes in $Rep(X-x,y, n, \tau)$. \\
  

$URep(X-x,y, n, \tau)$ & unitary representations in $Rep(X-x,y, n, \tau)$. \\


$UR^{irr}(X-x,y, n, \tau)$ & irreducible representations in $UR(X-x,y, n, \tau)$.\\


$UR^{irr}(X-x, n, \tau)$ &  conjugacy classes in  $UR^{irr}(X-x,y, n, \tau)$.\\

$Con(X-x, y, n, \tau)$ & isomorphism classes of triples $(F,\nabla, y^*)$ where \\
                       & $F$ is a holomorphic vector bundle on $X-x$ of rank $n$, \\
                       & $y^* : F_y \stackrel{\sim}{\to} \C^n$  a frame over $y$, \\
                       & and $\nabla : F \to \Omega_X^1\otimes F$ a holomorphic connection, \\
                    & with local monodromy $\rho(c) = \tau I$ where $\tau \in \C^{\times}$ and \\
                  &  $\rho : \pi_1(X-x,y) \to GL(n)$ denotes the monodromy \\
                  & representation of $(F, \nabla, y^*)$.\\

$Con(X-x, n, \tau)$ & isomorphism classes of pairs $(F,\nabla)$ where \\
                    & $F$ is a holomorphic vector bundle on $X-x$ of rank $n$, \\
                    & and $\nabla : F \to \Omega_X^1\otimes F$ a holomorphic connection, \\
                    & with local monodromy $\rho(c) = \tau I$ where $\tau \in \C^{\times}$.\\


$UCon(X-x, y, n, \tau)$ & connections in $Con(X-x, y, n, \tau)$ \\
                        & with unitary monodromy.\\

$UC^{irr}(X-x, y, n, \tau)$ & connections in $UCon(X-x, y, n, \tau)$ \\
                           & with irreducible monodromy.\\

$UC^{irr}(X-x, n, \tau)$ & the subset of $Con(X, x, n, \tau)$ where \\
                        & the monodromy (up to conjugacy) is irreducible and unitary.\\

$Log(X, x, y, n, \lambda)$ & isomorphism classes of triples $(E, \nabla, y^*)$ where\\
                           & $E$ is a holomorphic vector bundle on $X$ of rank $n$,\\
                         & $y^*$ is a frame for $E$ over $y$, \\
                       & and $\nabla : E \to \Omega_X^1(\log x)\otimes E$ a logarithmic connection \\
                   & with residue $\res_x(\nabla) = \lambda I \in End(E_x)$ where $\lambda \in \C$.\\

                             
$Log(X, x, n, \lambda)$ & isomorphism classes of pairs $(E, \nabla)$ where\\
                        &   $E$ is a holomorphic vector bundle on $X$ of rank $n$,\\
                  & and $\nabla : E \to \Omega_X^1(\log x)\otimes E$ a logarithmic connection \\
                & with residue $\res_x(\nabla) = \lambda I \in End(E_x)$  where $\lambda \in \C$.\\


$ULog(X, x, y, n, \lambda)$ & the subset of $Log(X, x, y, n, \lambda)$ where \\
                           & the monodromy is unitary. \\


$UL^{irr}(X, x, y, n, \lambda)$ & the subset of $ULog(X, x, y, n, \lambda)$ where \\
                              & the monodromy is irreducible.\\

$UL^{irr}(X, x, n, \lambda)$ & the subset of $Log(X, x, n, \lambda)$ where \\
           & the monodromy (up to conjugacy) is irreducible and unitary.\\ 
                             
\end{tabular}
                             
}

%
%

\stm\label{conjugacy} {\bf Remarks on conjugacy of representations.}
(1) If two unitary representations $\rho_1,\rho_2: \Gamma \to U(n)$ 
of a group $\Gamma$ are conjugate by an element $A\in GL(n)$, then it follows by complete
reducibility and Schur's lemma
that $\rho_1$ and $\rho_2$ are conjugate by an element $B\in U(n)$.
(2) Also by Schur's lemma, the inclusion $U(n)\hra GL(n)$
induces a bijection from 
the set of all conjugacy classes (under conjugacy by $U(n)$) of irreducible 
representation $\rho : \pi_1(X-x,y) \to U(n)$ with local monodromy $\rho(c) = \exp(-m/n)$
to 
the set of all conjugacy classes (under conjugacy by $GL(n)$) of irreducible 
representations $\rho : \pi_1(X-x,y) \to GL(n)$ with local monodromy $\rho(c) = \exp(-m/n)$,
such that $\rho$ is `unitarizable', that is, $\rho$ is conjugate under $GL(n)$ to a
unitary representation.\\
(3) If $y_0,y_1 \in X-x$, and $\gamma : [0,1]\to X-x$ is a path with $\gamma(0) = y_0$, 
$\gamma(1) = y_1$, then there is an induced isomorphism
$[\gamma]_* : \pi_1(X-x,y_0) \to  \pi_1(X-x,y_1)$ that depends only on the
path homotopy class of $[\gamma]$.
If $\delta : [0,1]\to X-x$ is another such path, then $[\delta]_*$ is the conjugate
of $[\gamma]_*$ by the element $[\alpha] = [\delta^{-1}\circ \gamma] \in \pi_1(X-x,y_0)$.
Hence the set of all conjugacy classes under $GL(n)$ of representations in
$Rep(X-x,y_0,n)$ is in a canonical bijection with the
set of all conjugacy classes under $GL(n)$ of representations in
$Rep(X-x,y_1,n)$, independently of the choice of $\gamma$. Hence we drop $y$ from the notation,
and denote the set of all conjugacy classes under $GL(n)$ in 
$Rep(X-x,y, n, \tau)$ simply by $Rep(X-x, n, \tau)$. Similar remark applies to the
other sets in the table above, where we have dropped $y$ from the notations.

\section{ODEs and the Deligne extension}

Let $U\subset \C$ be an open disk with center $0$. Let $x\in U$ be the point $z=0$. 
Any logarithmic connection $\nabla : \OO_U \to \Omega^1_U(\log x)\otimes \OO_U$
is uniquely determined by specifying  
  $$\nabla(e_j) = \frac{dz}{z}\otimes (\textstyle\sum_i \, A^i_j(z)e_i).$$  
where $A(z) = (A^i_j(z))$ is any $n\times n$-matrix
  of holomorphic functions on $U$, and $e_i$ is the standard free basis of $\OO_U^n$. 
  The residue of $\nabla$ at
  $z=0$ is given by the matrix $R = A(0)$. In classical terms,
  such a logarithmic connection corresponds to the linear system of ODEs
  $$\frac{d F(z)}{d z} = - \frac{A(z)}{z} \,F(z)$$
  where $F(z)$ is an $n\times 1$-vector of functions of $z$.

  \stm\label{exp of res is mono} {\bf Local monodromy for a good residue.} 
{\it With notation as above, 
    suppose that no two eigenvalues of $R = A(0)$ differ by a non-zero integer. 
    Then there exists a unique invertible $n\times n$-matrix of holomorphic functions
    $P(z)$ on $U$
with $P(0) = I$ such that $P(z)z^{-R}$ is a fundamental solution matrix for the above system, that is,
$$z \frac{d }{d z}(P(z)z^{-R}) = -A(z)P(z).$$
Consequently, the monodromy for the connection around $z=0$ is $e^{- 2\pi i R}$.
} 


\proof This is just Coddington-Levinson [C-L] Chapter 4, Theorem 4.1. Note that
the residue as defined in [C-L] is $-1$ times the residue $R$ as per our convention here.
Recall that in [C-L], the matrix valued function $z^{-R} = e^{-R \log z }$
is regarded as a multivalued holomorphic function of $z$, which becomes single valued over the 
universal cover of $U-x$ with coordinate $\log z$. Hence the monodromy is 
the multiplying matrix factor $e^{-2\pi i R}$ on the right by which the fundamental solution
$Pz^{-R}$ is modified when $\log z$ changes from $0$ to $2\pi i$ as $z$
travels once around $z=0$ in an anticlockwise sense.
\hfill$\square$

\stm\label{Deligne construction} {\bf Deligne construction of logarithmic extension.} 
Let $X$ be a Riemann surface (not necessarily compact), and let $x,y\in X$, with $x\ne y$.  
Let $(F,\nabla,y^*)$ be a framed holomorphic connection of rank $n$ over $(X-x,y)$. 
Let $\rho(c) = e^{-2\pi i R}\in GL(n)$ where $R\in M(n)$ is an $n\times n$ matrix over $\C$,
such that no two distinct eigenvalues of $R$ differ by a nonzero integer.
Then there exists an extension $(E, \phi: E|_{X-x} \stackrel{\sim}{\to} F)$ of $F$ to
a vector bundle $E$ on $X$, such that $\nabla$ extends as a logarithmic connection
on $E$ such that the residue $\res_x(\nabla) \in End(E_x)$ is conjugate to $R \in M(n)$.
Moreover, such an extension $(E, \phi)$ of $F$ is unique up to a unique isomorphism
of extensions. To construct $(E,\phi)$, consider the bundle $\OO^n_U$ on the disk $U\subset X$
chosen as in the beginning of Section 2. Let
$\nabla_U : \OO^n_U \to \Omega^1_U(\log x)\otimes \OO^n_U$ be the logarithmic connection defined by
$\nabla(e_i) = (dz/z) \otimes \sum_j R^j_ie_j$, where $e_i$ is the standard basis of $\OO^n_U$.
Note that $e_i$ defines a frame $(e_{i,y})$ over $y\in U$.
Then a fundamental solution of $\nabla_U$ on $U-x$ is given by the multivalued matrix
$z^{-R} = e^{- R\log z}$,
and it takes the values $I$ at $\log z =0$ and $e^{-2 \pi i R}$ at $\log z = 2\pi i$.
This shows that for the monodromy $\rho_U(c) = e^{-2\pi i R}$. Hence there exists
a unique isomorphism of framed holomorphic connections
$$\psi : (F, \nabla, y^*)|_{U-x} \to (\OO_U, \nabla, (e_{i,y}))|_{U-x}$$
over $U-x$.
Gluing $F$ and $\OO_U$ over $U-x$ via $\psi$ defines the extension $(E, \phi)$.
The uniqueness of $(E, \phi)$ follows from Lemma \ref{exp of res is mono} applied to
the logarithmic connection $(\un{Hom}(E_1,E_2), \nabla)$
made from two extensions $(E_1,\nabla_1)$ and $(E_2,\nabla_2)$,
which has the property that no nonzero integer is an eigenvalue of $\res_x(\nabla)$,
as by assumption, no two eigenvalues of $R$ differ by a nonzero integer.
This has the consequence that the integrable section $\sigma$ of $\un{Hom}(E_1,E_2)$
on $X-x$ defined by $\id_F$ extends to $x$ to define a flat isomorphism $E_1\to E_2$.
(See [3] or Theorem 4.4 in [6] for more details).

\stm {\bf Note.} We will write the extension $(E,\phi)$ simply as $E$, and we will say that
`$E$ is the logarithmic extension of $(F, \nabla)$ with residue conjugate to $R$', suppressing $\phi$
from the notation unless it is explicitly needed for unambiguity. 

\stm {\bf Example.} It is necessary in the above to assume that no two eigenvalues of $R$
differ by a nonzero integer. For example, let $F = \OO^2_{X-x}$, let $\nabla(e_i) =0$
where $e_1,e_2$ form the standard basis of $\OO^2_{X-x}$, and let $y^*$ be the basis 
$(e_{i,y})$ over $y$. Let $R = \diag(1,0)$. Then $F$ has two possible non-isomorphic extensions
$\I_x\oplus \OO_X$ and $\OO_X \oplus I_x$, on which $\nabla$ extends logarithmically with 
residues conjugate to $R$. (Note the isomorphic bundles $\I_x\oplus \OO_X$ and $\OO_X \oplus I_x$ are
not isomorphic as {\it extensions} of $F$.)


\stm\label{Deligne for scalar local monodromy} {\bf Case of scalar local monodromy.} 
Let $X$ be a Riemann surface (not necessarily compact), and let $x,y\in X$, with $x\ne y$.  
Let $(E,\nabla,y^*) \in Log(X,x,y,n, \lambda)$, that is, 
$E$ is rank $n$ bundle with a logarithmic connection
$\nabla : E \to \Omega^1_X(\log x)\otimes E$ and a frame $y^*$ for $E$ over $y$,
with residue $\res_x(\nabla) = \lambda I \in End(E_x)$ for a given
$\lambda \in \C$. 
Then its restriction $(E|_{X-x}, \nabla|_{X-x}, y^*)$ to $X-x$ lies in 
$Con(X-x,y,n, e^{-2\pi i \lambda})$, that is, it is a framed
holomorphic connection on $(X-x,y)$, whose local monodromy around $x$ is
$$\rho(c) = e^{-2 \pi i \lambda}I \in GL(E_y) \stackrel{y^*}{=} GL(n)$$
Hence restriction defines a map
$${\cal R}_y : Log(X,x,y,n,\lambda) \to Con(X-x,y,n, e^{-2\pi i \lambda})$$
that sends $(E,\nabla, y^*)$ to $(E|_{X-x}, \nabla|_{X-x}, y^*)$.
The Deligne construction \ref{Deligne construction} gives its inverse map
$$\D_y : Con(X-x,y,n, e^{-2\pi i \lambda})\to Log(X,x,y,n,\lambda)$$
showing that ${\cal R}_y$ and $\D_y$ are bijections. 
Forgetting the point $y$ and the frame $y^*$, the above map descends to a bijective map
${\cal R} : Log(X,x,n,\lambda) \to Con(X-x,n, e^{-2\pi i \lambda})$
with inverse bijection
$\D : Con(X-x,n, e^{-2\pi i \lambda})\to Log(X,x,n,\lambda)$.

\section{The case of line bundles}

\stm\label{fun grp of pun rie sur}
{\bf The fundamental group of a punctured Riemann surface.}
As before, let $X$ be a compact Riemann surface of genus $g\ge 0$, and let 
$x,y\in X$ with $x\ne y$. Let there be given 
a local holomorphic coordinate chart $(U,z)$ on $X$ such that
the image $z(U)\subset \C$ is an open disk of radius $>1$,
and let $x,y\in U$ with $z(x)=0$, $z(y) = 1$.  
The element $c\in \pi_1(X-x, y)$ is then defined by the loop $z(t) = e^{2 \pi i t}$ in $U - x$.
It is well known that there exist $2g$ distinct elements
$a_1,\ldots, a_g, b_1,\ldots, b_g \in \pi_1(X-x, y)$ such that
$\pi_1(X-x, y)$ is the free group on these $2g$ elements,
and the local positive loop $c$ is given by 
$$ c = a_1b_1a_1^{-1}b_1^{-1} \cdots a_gb_ga_g^{-1}b_g^{-1}.$$
As $c$ is a product of commutators, its homology class is zero, that is, 
$[c] = 0\in H_1(X-x) = \pi_1(X-x, y)/(\pi_1(X-x, y),\pi_1(X-x, y))$.
Note that the inclusion $X-x \hra X$
induces an isomorphism of groups
$$\pi_1(X-x, y)/\langle c \rangle \stackrel{\sim}{\to} \pi_1(X, y)$$
where $\langle c \rangle$ denotes the smallest normal subgroup of $\pi_1(X-x, y)$
that contains $c$.

For proving the Main Theorem \ref{main theorem}, 
we begin with line bundles on $X$, that is, we take $n=1$. 
Note that every line bundle is a stable vector bundle.

As the group $U(1)$ is abelian and as $c$ lies in the commutator subgroup of 
$\pi_1(X-x,y)$, for any $\rho : \pi_1(X-x,y) \to U(1)$ we must have $\rho(c)=1 \in U(1)$.
Hence $\rho$ uniquely factors
via the quotient $\pi_1(X-x,y) \to \pi_1(X,y)$, to define
a representation which we again denote by $\rho : \pi_1(X,y) \to U(1)$.
In fact, as $U(1)$ is abelian, the choice of the base point $y$ does not
matter, and $\rho$ further factors uniquely via the quotient $\pi_1(X,y) \to H_1(X)$
to define a representation,  
again denoted by $\rho : H_1(X)\to U(1)$.
%
%
%
%
Also, for any holomorphic connection $\nabla$ on a line bundle $L$ on $X$,
the choice of a base point $y\in X$
and a basis $y^*$ of the fiber $L_y$ over it does not matter
for the corresponding monodromy representation
$\rho : H_1(X)\to U(1)$, which remains independent of such choices.

\stm \label{all line bundles of deg zero are uniquely unitary}
{\bf Theorem} {\it Let $X$ be a compact Riemann surface.
Then for any holomorphic line bundle $L$ on $X$ with $\deg(L)=0$, 
there exists a unique representation 
$\rho : \pi_1(X) \to U(1)$ such that $L$ is associated to the representation $\rho$.
Equivalently, there exists a unique holomorphic connection $\nabla$
on $L$ such that the monodromy of $\nabla$ is unitary.}

{\bf Proof} (following Narasimhan's Trieste lecture notes [9]).
We have the following commutative diagram of sheaves on $X$, with exact rows,
where for any of the groups group $A = \Z,\, \R,\, U(1)$, we denote by $A_X$ the sheaf on $X$ 
of germs of locally constant functions with values in $A$ (which effectively means that we have
given the various groups $A$ the discrete topology).

$$%
\begin{array}{ccccc}
0\to & \Z_X & \to \R_X \to & U(1)_X & \to 0 \\
     & \|   & \da          & \da    & \\
0\to & \Z_X & \to \OO_X \to & \OO_X^{\times} &  \to 0
\end{array}$$

As the boundary map $\pa: Pic(X) = H^1(X,\OO_X^{\times})\to H^2(X,\Z_X) = \Z$
sends \\
$[L]\mapsto \deg(L)$, we obtain another commutative diagram with exact rows

$$
\begin{array}{ccccc}
H^1(X, \Z_X) & \to H^1(X, \R_X) \to & H^1(X, U(1)_X) & \to 0 \\
    \|       & \da                  & \da            &       \\
H^1(X, \Z_X) & \to H^1(X, \OO_X)\to & Pic^0(X)       & \to 0 
\end{array}$$

The map $H^1(X,\R_X) \to H^1(X,\OO_X)$ is an isomorphism as $X$ is compact K\"ahler, 
hence $H^1(X,U(1)_X)  \stackrel{\sim}{\to} Pic^0(X)$.
By the universal coefficient theorem, the natural map 
$Hom(\pi_1(X), U(1)) = Hom(H_1(X), U(1)) \to 
H^1(X,U(1)_X)$ is an isomorphism. Hence, as claimed, the map 
$Hom(\pi_1(X), U(1))\to Pic^0(X)$ is an isomorphism.
\hfill$\square$

\stm {\bf Example.} The unique connection on $L = \OO_X$ with unitary monodromy
is the de Rham connection (exterior derivative) $d : \OO_X \to \Omega^1_X$.

\stm\label{log connection on O(m x)}
{\bf The natural logarithmic connection on the line bundle $\OO_X( m x)$.}
Note that $\Omega^1_X(\log x)$ is naturally isomorphic to $\Omega^1_X\otimes_{\OO_X} \OO_X(x)$.
The de Rham connection $d : \OO_X \to \Omega^1_X$ maps the ideal sheaf
$\I_x \subset \OO_X$ 
into $\Omega_X = \Omega_X(\log x)\otimes_{\OO_X} \I_x$, which is 
a logarithmic connection on $\I_x = \OO_X(-x)$. Taking duals and tensor powers, this induces 
a logarithmic connection $\nabla_m$ on the line bundle $\OO_X( m x)$ for any $m\in \Z$.
These logarithmic connections have the following coordinate description
in terms of any local coordinate chart $(U,z)$ around $x$ in which $x$ is given by $z=0$.
Note that $\OO_X(mx)|_{X-x} = \OO_{X-x}$, and $\nabla_m|_{X-x} = d$, the de Rham
connection on $\OO_{X-x}$. Over $U$, the sheaf $\OO_X(mx)$ has the free basis $z^{-m}$, and
a simple calculation shows that
$$\nabla_m (z^{-m}) = - m \frac{dz}{z}\otimes z^{-m}.$$
In particular, the residue of $\nabla_m$ at $x$ is $-m$, and so we have
$\deg(\OO_X(m x)) = m = - \trace \res_x(\nabla_m)$,
as given by Lemma \ref{deg is minus res}.
As $\nabla_m$ restricts to $d$ on $X-x$, the constant functions on $X-x$ are
flat sections of $\OO_X(m x)$ over $X-x$, and so 
the monodromy of $\nabla_m$ is trivial (that is, $\rho$ takes the constant value $1$).
In particular, the monodromy of $\nabla_m$ is unitary.

\stm\label{log con on line bundles}
{\bf Arbitrary line bundles of degree $m$.} If $\deg(L) = m$,
then $\deg(L\otimes \OO_X(-mx)) = 0$, so by
Theorem \ref{all line bundles of deg zero are uniquely unitary},
there exists a unique holomorphic connection $\nabla_{L\otimes \OO_X(-mx)}$ on it
whose monodromy is unitary. Tensoring the holomorphic connection
$(L\otimes \OO_X(-mx), \nabla_{L\otimes \OO_X(-mx)})$ with 
$(\OO_X(mx), \nabla_m)$, we see that $L$ has a unique logarithmic connection on $(X,x)$ with unitary
monodromy, which has trivial local monodromy $\rho(c) = e^{2\pi i m/1} = 1$ and residue
$-m = -\mu(L)$. 

This concludes the proof of Theorem \ref{main theorem} for rank $1$ and arbitrary degree $m$.

\stm {\bf Relation between residue and degree in all ranks.}

The following lemma, and its appropriate generalization to varieties of higher dimensions, are
well known. For easy reference, we include a proof for compact Riemann surfaces.

\stm\label{deg is minus res} {\bf Lemma.} {\it
Let $X$ be a compact Riemann surface, and let $x\in X$.
If $(E,\nabla)$ is a logarithmic connection on $(X,x)$,
then $\deg(E) = - \res_x(\nabla)$.
} 

{\bf Proof.}  
The cohomology sequence of the
exact sequence $0\to \Omega^1_X\to \Omega^1_X(\log x) \to \OO_x \to 0$, together
with the fact that $H^1(X, \Omega^1_X(\log x)) =0$ by Serre duality, shows that 
the induced map $H^0(X,\Omega^1_X) \to H^0(X,\Omega^1_X(\log x))$ is an isomorphism.
The logarithmic connection
$\nabla : E \to \Omega^1_X(\log x)\otimes E$ induces
a logarithmic connection $\nabla' : \det(E) \to \Omega^1_X(\log x)\otimes \det(E)$.
Note that
$$\res_x(\det (E), \nabla') = \trace \,\res_x(E,\nabla).$$
If $\deg(E) =m$, then the line bundle
$L = \det (E) \otimes \OO_X(-m x)$ has a logarithmic connection $\nabla''$ 
induced by the logarithmic connections on $\det (E)$ and $\OO_X(-m x)$.
Note that 
$$\res_x(L, \nabla'') = \res_x (\det (E), \nabla') + 
\res_x (\OO_X(-m x), \nabla_{- m}) = \trace \,\res_x(E,\nabla) - m.$$
As $L$ is a line bundle of degree $0$,
by Lemma \ref{all line bundles of deg zero are uniquely unitary} 
it admits a holomorphic connection $D : L \to \Omega^1_X\otimes L$.
Hence $\nabla '' - D \in H^0(X, \Omega^1_X(\log x)) = H^0(X, \Omega^1_X)$. 
This shows that $\res_x(L, \nabla '') = \res_x(L, D) = 0$.
It follows that $\trace \,\res_x(E,\nabla) = m$
\hfill$\square$

\section{Logarithmic reformulation of Weil's theorem}

The statement  (2) of the following theorem is a
generalization of Weil's theorem to all degrees, which in another equivalent form
already occurs as Proposition 6.2
in [13]. We give a (more transparent) logarithmic reformulation of both the statement and the proof.
I had reported the statement (1) in a talk in Strasbourg in 1992, but did not publish it. 


\stm\label{generalized Weil} {\bf Theorem.} 
{\it Let $X$ be a compact Riemann surface, and let $x\in X$.
  Then we have the following.

{\bf (1)} Every holomorphic vector bundle $E$ on $X$ admits a logarithmic connection \\
$\nabla : E  \to \Omega_X^1(\log x)\otimes E$.

{\bf (2)} For any holomorphic vector bundle $E$ on $X$, the following three statements are equivalent:\\
$~~~~$ (a) $E$ admits a logarithmic connection $\nabla$ with residue $\res_x(\nabla) = - \mu(E) I$.\\  
$~~~~$ (b) For every direct summand $F$ of $E$ we have $\mu(F) = \mu(E)$.\\ 
$~~~~$ (c) The slope of every indecomposable component of $E$ equals the slope of $E$.
} 

\proof This is a logarithmic version of Atiyah's proof of Weil's theorem, with
some modifications. We assume familiarity with Atiyah's proof (see [1]). 
The obstruction for the existence of a holomorphic
connection on $E$ is the Atiyah class $At(E) \in H^1(X,\Omega^1_X\otimes \un{End}(E))$. 
If $E$ is described by an open cover $(U_i)$ and transition functions
$(g_{ij})$, then $At(E)$ is induced by the Cech $1$-cocycle
$(dg_{ij}g^{-1}_{ij})$ w.r.t. this cover. The analogous argument for logarithmic
connections shows that the obstruction for the existence of a logarithmic
connection on $E$ over $(X,x)$ is (what we will call as) the {\it logarithmic} Atiyah class
$\alpha_E \in H^1(X,\Omega^1_X(\log x)\otimes \un{End}(E))$, described by 
the Cech $1$-cocycle $(dg_{ij}g^{-1}_{ij})$.

To prove (1), we will show that $\alpha_E =0$ for each $E$. 
As $\un{End}(E)$ is self-dual under the trace pairing, the pairing
$$H^0(X, \OO_X(-x)\otimes \un{End}(E)) \times H^1(X,\Omega^1_X(\log x)\otimes \un{End}(E))
\stackrel{trace}{\to} H^1(X,\Omega^1_X) = \C$$
is non degenerate by Serre duality. Note that $H^0(X, \OO_X(-x)\otimes \un{End}(E))$
is the subspace of $H^0(X, \un{End}(E)) = End(E)$ consisting of all
endomorphisms $\varphi : E\to E$ such that 
$\varphi$ vanishes at the point $x$. As $X$ is compact, the characteristic polynomial
$P(\varphi_u) \in \C[T]$ 
of $\varphi_u = \varphi|_{E_u}$ on any fiber $E_u$ is independent of $u \in X$, and hence 
$P(\varphi_u)= P\varphi_x) = T^n$, where $n= \rank(E)$. Hence $\varphi_u^n=0$ for all $u$, so
$\varphi^n=0$. Hence as in [1], there exists
a nested sequence (full flag)
of vector subbundles $0 = E_0\subset E_1\subset\ldots \subset E_n=E$,
such that $\varphi(E_i) \subset E_{i-1}$ for all $1\le i\le n$.
We can choose an open cover $U_i$ of $X$ with trivializations
of $E$ over $U_i$, such that the transition functions $g_{ij}$ preserve this full flag.
Hence $dg_{ij}g_{ij}$ are upper triangular matrices while $\varphi$ is given by strictly
upper triangular matrices $\varphi_i$ over $U_i$. Hence $\varphi_i dg_{ij}g^{-1}_{ij}$
is again strictly upper triangular, so its trace is zero. 
This shows that $\langle\varphi,\,\alpha_E \rangle = 0$ for all
$\varphi\in H^0(X, \F \otimes \I_x)$, hence $\alpha_E=0$. Hence (1) is true. 

Next we prove (2). A logarithmic connection on $E = E'\oplus E'$ over $(X,x)$ with residue $\lambda I$ 
induces by inclusions and projections logarithmic connections on $E'$ and $E''$ over $(X,x)$
with residues $\lambda I$ (here, the same notation $I$ stands for the identity endomorphisms
of $E_x$, $E'_x$, $E''_x$ for simplicity). Also, by Krull-Remak theorem, 
$E$ is the direct sum of indecomposable subbundles. So to prove (2), it is enough to show
that if $E$ is indecomposable then $E$ admits a logarithmic connection $\nabla$ over $(X,x)$ with
$\res_x(\nabla) = \lambda I$ for some $\lambda\in \C$.  Then we will necessarily  
have $\lambda= - \mu(E)$ by Lemma \ref{deg is minus res}.

Let $ev_x : \un{End}(E) \to (\OO_X/\I_x)\otimes_{\OO_X} \un{End}(E)  = End_{\C}(E_x)$
be the evaluation map, where $\I_x\subset \OO_X$ denotes 
the ideal sheaf of $x$ in $X$. Let $\C I \subset End_{\C}(E_x)$ be the scalar
%
%
endomorphisms.
Let $\F \subset \un{End}(E)$ be the coherent sub-$\OO_X$-module that is the
inverse image of $\C I$ under $ev_x$. Note that $\F$ is locally free, but it is
not a subbundle of $\un{End}(E)$.

Repeating the arguments of [1] in this context shows that 
the obstruction to the existence of a logarithmic connection on $E$ over $(X,x)$ with
residue of the form $\lambda I$ is defined by the $1$-cocycle
$\beta_E = (dg_{ij}g^{-1}_{ij})$ in $H^1(X, \Omega^1_X(\log x) \otimes \F)$.

The sheaf $\F$ fits in the exact sequence
$$0\to \un{End}(E)\otimes \I_x \to \F \to \OO_x \to 0$$
where $\OO_x = \OO_X/\I_x$, and the map $\F \to \OO_x$ sends an endomorphism $\varphi$
to $\lambda$ where $\varphi_x = \lambda I$. This means
the above map $\F \to \OO_x$ sends $\varphi$ to $(1/n)\trace (\varphi_x)$. 
Tensoring with $\OO_X(x)$ gives the exact sequence
$$0\to \un{End}(E) \to \F\otimes \OO_X(x) \to \OO_X(x)_x \to 0.$$
As $\un{End}(E)$ is self-dual under the trace pairing, applying the functor
$\un{Hom}(-,\OO_X)$ 
we get an exact sequence of sheaves 
$$0\to \F^{\vee}\otimes \I_x \to \un{End}(E) \to \OO_x\to 0$$
where $\un{Ext}^1(\OO_X(x)_x, \OO_X)$ is canonically identified with $\OO_x$,
and the map $\un{End}(E) \to \OO_x$ sends $\varphi$ to $(1/n)\trace (\varphi_x)$. 
The long exact cohomology sequence of the above short exact sequence of sheaves shows that 
$$H^0(X, \F^{\vee} \otimes \I_x) = \{ \varphi \in End(E)\,|\, \trace(\varphi_x) = 0 \}.$$

By Serre duality we have a nondegenerate pairing
$$\langle\mhyphen,\,\mhyphen\rangle :
H^0(X, \F^{\vee}\otimes \I_x) \times H^1(X, \Omega^1(\log x)\otimes \F) \to H^1(X, \Omega^1) =\C,$$
therefore to show that $\beta_E = 0 \in H^1(X, \Omega^1(\log x)\otimes \F)$,
we must show that $\langle\varphi, \beta_E\rangle = 0$ for all
$\varphi \in End(E)$ such that $\trace(\varphi_x) = 0$.
In cocycle terms, we must show that for any such $\varphi$, 
$\trace(\varphi_i \circ dg_{ij}g_{ij}^{-1}) = 0 \in H^1(X, \Omega^1)$.

As by assumption $E$ is indecomposable, each global endomorphism $\varphi$ of $E$
is of the form $c I + N$ where $c\in \C$ and $N$ is nilpotent.
If $\varphi = cI$, then requirement that $\trace(\varphi_x)=0$ means $cI=0$.
Hence $\varphi$ is nilpotent. Then as above, we can choose a nested sequence (full flag)
of subbundles $0= E_0\subset E_1\subset\ldots \subset E_n=E$, where $n= \rank(E)$,
such that $\varphi(E_i) \subset E_{i-1}$ for all $1\le i\le n$.
We can choose an open cover $U_i$ of $X$ with trivializations
of $E$ over $U_i$, such that the transition functions $g_{ij}$ preserve this full flag.
Hence $dg_{ij}g_{ij}$ are upper triangular matrices while $\varphi$ is given by strictly
upper triangular matrices $\varphi_i$ over $U_i$. Hence $\varphi_i dg_{ij}g^{-1}_{ij}$
is again strictly upper triangular, so its trace is zero. 
This shows that $\langle\varphi,\,\beta_E \rangle = 0$ for all
$\varphi\in H^0(X, \F \otimes \I_x)$, hence $\beta_E=0$. \hfill$\square$

\bigskip

Combining with the above theorem with the bijection
of statement \ref{Deligne for scalar local monodromy} gives us the following.


\stm\label{generalized Weil again}
{\bf Theorem.} {\it Let $X$ be a compact Riemann surface, and let $x\in X$.
Then we have the following.

{\bf (1)} For any indecomposable vector bundle $E$ on $X$, there exists an indecomposable representation 
$\rho : \pi_1(X -x) \to GL(n)$ with local monodromy $\rho(c) = e^{- 2\pi i \mu(E)}I$ around $x$,
such that $E$ is the underlying vector bundle of the logarithmic connection
$(E,\nabla)$ on $(X,x)$ associated to $\rho$ by the Deligne construction made by choosing
$- \mu(E) I$ as the residue of the logarithmic connection at $x$.

{\bf (2)}  Conversely, if $\rho : \pi_1(X -x) \to GL(n)$ is indecomposable, with
$\rho(c) = e^{2\pi i m/n}I$ where $m, n\in \Z$ with $n\ge 1$, then the vector bundle $E$ on $X$ which is 
given by the Deligne construction made by choosing
$(- m/n) I$ as the residue of the logarithmic connection at $x$,  
is an indecomposable vector bundle on $X$ of rank $n$ and degree $m$.
}   %

\stm {\bf Existence of indecomposable bundles.}
If $g(X) = 0$ then $\pi_1(X-x,y)$ is trivial,
the only indecomposable representation of $\pi_1(X-x,y)$ is the trivial representation
of rank $1$, and the only indecomposable vector bundles
on $X$ are of rank $1$ by the theorem of Grothendieck. 

When $g(X) = 1$, $\pi_1(X-x,y)$ is the free group on $2$ generators.
As proved by Atiyah, when $g(X)=1$ and $n\ge 1$, 
there exists an indecomposable vector bundle $E$ on $X$ of rank $n$ and any prescribed degree $m$.
This implies by the above Theorem \ref{generalized Weil again}
that there exists an indecomposable rank $n$ representation of $F_2$
with $\rho(c) = e^{2\pi i m/n}I$, for any $m$, $n$ with $n\ge 1$. 
Venkataramana provided a direct example of such a representation
(see Lemma \ref{existence of irr reps genus one} below). This, combined with
Theorem \ref{generalized Weil again}, directly shows that indecomposable bundles exist
on elliptic curves for any given rank $n\ge 1$ and degree $m$, bypassing the results
of $[2]$.

When $g(X)\ge 2$, [13] Proposition 9.1 constructs 
irreducible rank $n$ unitary representations of the free group $\pi_1(X-x,y)$ with 
$\rho(c) = e^{2\pi i m/n}I$, for any value of $m$ and $n$ with $n\ge1$.
They then deduce the existence of stable bundles
of rank $n$ degree $m$ from it by one part of the Narasimhan-Seshadri theorem. This part
is then used to prove the remaining part of the theorem. Of course, all stable bundles
are indecomposable.

\section{Logarithmic reformulation of the Narasimhan-\-Seshadri theorem}

Any stable bundle $E$ on a compact Riemann surface $X$
is indecomposable. So by Theorem \ref{generalized Weil again},
$E$ arises from some indecomposable representation $\rho : \pi_1(X-x,y) \to GL(n)$
with scalar local monodromy by Deligne construction 
with scalar residue. However, is there some condition on the monodromy representation,
stronger than indecomposability, 
that guarantees that the bundle $E$ is stable? Does every stable bundle arise
this way from a monodromy representation that satisfies this stronger condition?
This (and more) was historically answered by the Narasimhan-Seshadri theorem.

The following, which is the main result of this note, is an equivalent logarithmic
%
%
reformulation of the Narasimhan-Seshadri theorem. We allow any value of $g$.

\stm\label{main theorem} {\bf Main Theorem.} {\it
Let $X$ be a compact Riemann surface of any genus $g\ge 0$, and let $x,y\in X$ with $x\ne y$. Let 
$n, m\in \Z$ with $n\ge 1$. Then we have the following.

    (1) If $\rho : \pi_1(X-x,y) \to U(n)$ is an irreducible unitary representation
  whose local monodromy around the puncture $x$ equals $e^{2\pi i m/n}I$, then the  
  Deligne extension $E(\rho)$ of $E_{\rho}$,
  with residue $(-m/n)I$ at $x$,
  is a stable holomorphic vector bundle on $X$ of rank $n$ and degree $m$. Moreover,
  two such representations $\rho_1$ and $\rho_2$ are conjugate if and only if the bundles 
  $E(\rho_1)$ and $E(\rho_2)$ are isomorphic.
 
%

  (2) If $E$ is a stable holomorphic vector bundle on $X$ of rank $n$ and degree $m$,
  then there exists an irreducible representation $\rho : \pi_1(X-x,y) \to U(n)$, 
  with local monodromy $e^{2\pi i m/n}I$ around $x$,   
  such that $E$ is isomorphic to $E(\rho)$. 

  (3) The resulting map $UR^{irr}(X-x, n, e^{2\pi i m/n}) \to M^s(X,n,m) : \rho\mapsto E(\rho)$
  is a real analytic isomorphism between  \\
$~~~$    (a) the manifold $UR^{irr}(X-x, n, e^{2\pi i m/n})$
  of all conjugacy classes of irreducible \\
$~~~$     unitary representations 
  $\rho : \pi_1(X-x) \to U(n)$ with local monodromy $e^{2\pi i m/n}I$ \\
$~~~$  around $x$, and\\
$~~~$    (b) the manifold $M^s(X,n,m)$ of all isomorphism classes of stable holomorphic \\
$~~~$ vector bundles $E$ on $X$ of rank $n$ and degree $m$.

} 

\stm The statement (1) is proved in Section 9. The strategy of proof is copied from [13],
but here the logarithmic reformulation leads to somewhat simpler arguments.
The statements (2) and (3) are almost identical to the corresponding statements
in [13], and we do not repeat the very beautiful proofs given there. 
In particular, we know by the arguments in [13] that $UR^{irr}(X-x, n, e^{2\pi i m/n})$ and
$M^s(X,n,m)$ are real analytic manifolds, and the bijective map 
$$UR^{irr}(X-x, n, e^{2\pi i m/n}) \to M^s(X,n,m)$$
of (3) is a tangent-level injection and hence an isomorphism of real analytic manifolds.
The original argument in [13] becomes somewhat simplified because of developments between
1965 and 1970, and this partially simplified proof of (2) and (3)
(which is standard knowledge since 1970s) is sketched in Section 10 for completeness.

\section{Comparison with Grothendieck's construction}

\stm\label{ass bun as inv dir im}{\bf Associated bundles as invariant direct images.}
Let $p: P\to X$ be a principal bundle under the left action of a group $\pi$, and
let $\sigma: \pi \to GL(n)$ be
a representation. The total space $P(\sigma)$ of the associated vector bundle $P(\sigma)\to X$
is the quotient $\pi \backslash (P\times \C^n)$, where the left action
$\pi \times (P\times \C^n) \to P\times \C^n$ is given by $g (y,v) = (gy, \rho(g)v)$.
The sheaf of sections $\un{P(\sigma)}$ of $P(\sigma)\to X$ is then the invariant direct image
$$\un{P(\sigma)} = p^{\pi}_*(\OO_P^n(\sigma))$$
where $\OO_P^n(\sigma)$ is the sheaf $\OO_P^n$ on $P$ with the obvious left $\pi$-action coming from
$\sigma$, that lifts the left $\pi$-action on $P$.


\stm\label{left right correspondence}{\bf Left-right correspondence.} Let $G = \pi^{op}$ be
the opposite group. Recall that $G$ has the same underlying set that $\pi$ has, but
the multiplication $*$ in $G$ is defined in the reverse order to that in $\pi$,
that is, $g* h = hg$. Any left $\pi$-action
$\pi \times Z \to Z : (g,z)\mapsto gz$ on a space $Z$ 
becomes a right $G$-action $Z\times G \to Z$ defined by $z* g = gz$. 
A left principal $\pi$-bundle $p: P\to X$ becomes a right principal $G$-bundle
$P \to X$ under the above action. A representation $\sigma : \pi \to GL(n)$ corresponds to
a representation $\rho = \sigma^{op} : G \to GL(n)$, defined by $\rho(g) = \sigma(g^{-1})$.
The vector bundle $P(\rho)$ is the quotient of $P\times \C^n$ by the
right $G$-action given by
$(y, v)* g = (y* g, \rho(g^{-1})v) = (gy, \sigma(g)v) = g(y,v)$. This shows that
$$P(\sigma) = P(\rho).$$

The formulas in \ref{ass bun as inv dir im} and \ref{left right correspondence}
will be used in our comparison between the Grothendieck construction of invariant direct images
and the Deligne construction of logarithmic extensions.

In [4], Grothendieck associates a bundle $\EE(\sigma)$ on a compact Riemann surface $X$ to
representations of a certain group $\pi$. 
Grothendieck's construction is based on the following lemma, which as
he says is essentially a topological fact. It is provable just from
the basics of fundamental groups and covering spaces, by attaching open disks
to a suitably constructed covering of $X-x$. 
It is possible that Poincar\'e proved the lemma via tessallations
of the hyperbolic plane, which he had invented.

\stm\label{Poincare ramified cover} {\bf Lemma.} (Poincar\'e.)
{\it Let $X$ be a compact Riemann surface of genus $g\ge 1$,
let $x\in X$, and let $N\ge 1$. Then there exists a Riemann surface
$Y$ with a holomorphic morphism $p: Y \to X$ which satisfies the following
properties.\\
(1) The space $Y$ is simply connected. \\
(2) The automorphism group $\pi = Aut(Y/X)$ acts properly discontinuously on the left on $Y$,
    and $p: Y \to X$ is the corresponding quotient.
    In particular, $\pi$ acts transitively on all fibers.\\
(3) The restriction $p': Y' = Y - p^{-1}(x) \to X -x = X'$ of $p$ is unramified. \\
(4) The inertia subgroup $\pi_{\wh{x}}$ at any point $\wh{x}\in p^{-1}(x)$
    is cyclic of order $N$. So if $N=1$ then 
    $p: Y\to X$ is a universal covering of $X$. 

Any $p_1: Y_1\to X$ that satisfies (1)-(4) above is isomorphic
to $p: Y\to X$, where the isomorphism is unique up to composition with an element of
the Galois group $\pi$.    
} 

\stm {\bf The requirement $g\ge 1$.} The above lemma is false when $g=1$ and $N\ge 2$.
It is to be noted that our reformulation of Narasimhan-Seshadri theorem via logarithmic connections
does not need
this lemma, and works uniformly for all values $g\ge 0$ of the genus.

\stm {\bf Relation with $\pi_1(X-x,y)$.} If $p: Y\to X$ is as above, then 
its restriction $p' : Y' = Y-p^{-1}(x) \to X -x = X'$ is an unramified Galois covering with
Galois group $Aut(Y/X)$. Let $y'\in Y'$ be a point over $y\in X'$.
Let $K\subset \pi_1(X',y)$ denote the image of $p'_* : \pi_1(Y',y')\to \pi_1(X',y)$.
Then $K$ is the smallest normal subgroup that contains $c^N \in \pi_1(X',y)$.
Moreover,
$$\pi^{op} = Aut(Y/X)^{op} = \pi_1(X',y)/K$$
where the opposite group occurs by our convention about composition of loops in fundamental groups.
As underlying sets of
a group and its opposite group are identical, 
the image $\ov{c} \in \pi_1(X',y)/K$ of $c$ defines an element $\ov{c}$ of $\pi = Aut(Y/X)$.


\stm\label{pi rep inn terms of fun grp rep}
{\bf Relationship between the representations of the groups $\pi_1(X-x,y)$ and $\pi$.} 
Any representation $\sigma : \pi \to GL(n)$ defines by composition 
a representation
$$\rho : \pi_1(X',y) \to  \pi_1(X',y)/K = \pi^{op} \stackrel{\sigma^{op}}{\to} GL(n)$$
where $\sigma^{op} : \pi^{op}\to GL(n)$ is defined by $\sigma^{op}(a) = \sigma(a)^{-1}$.
Conversely, any $\rho : \pi_1(X',y) \to GL(n)$ for which $\rho(c^N) = I$ arises
from a unique $\sigma:  \pi \to GL(n)$. 

Next suppose that $n\ge 1$ divides $mN$ where $m\in \Z$.
Then any representation $\rho : \pi_1(X-x,y) \to GL(n)$ with
$\rho(c) = e^{2\pi i m/n}I$ sends $c^N\mapsto I$, hence factors to define a representation
$\ov{\rho} : \pi_1(X-x,y)/K \to GL(n)$. We denote by $\ov{\rho}^{op} : \pi \to GL(n)$
the representation well-defined by the formula   
$$\ov{\rho}^{op}(\ov{\alpha}) = \rho(\alpha^{-1})\in GL(n)$$
where $\alpha\mapsto \ov{\alpha}$ under the quotient map $\pi_1(X-x,y) \to \pi_1(X-x,y)/K$.
Thus, every representation $\sigma: \pi  \to GL(n)$ for which
$\sigma(\ov{c}) = e^{2\pi i m/n}I$ is of the form 
$\sigma = \ov{\rho}^{op}$ for a unique $\rho: \pi_1(X-x,y) \to GL(n)$ with $\rho(c) = e^{2\pi i m/n}I$.

\stm {\bf Local description.} 
Let $\gamma : [0,1]\to X$ be an injective path with $\gamma(0) =y$, and $\gamma(1) = x$.
Then $\gamma$ has a unique lift $\wh{\gamma} : [0,1] \to Y$ with $\wh{\gamma}(0) = y'$.
Let $x' = \wh{\gamma}(1)$, which is a point in $p^{-1}(x)$. Then
the inertia subgroup $\pi_{x'}\subset \pi$ is the cyclic group 
$\{ e, \ov{c},\ldots, \ov{c}^{N-1}\}$.
The inverse image under $p$ of the open disk $U$ in $X$ (that was chosen at the beginning of Section 2)
is a disjoint union of open disks in $Y$ indexed by $p^{-1}(x)$, whose centers
are the points of $p^{-1}(x)$.
Let $W$ be the component of $p^{-1}(U)$ that contains $x'$. Then $W$ has a complex coordinate $w$
such that the map $p|_W : W\to U$ is given by $w\mapsto z = w^N$, and the right action of
$\ov{c} \in \pi_1(X',y)/K$ on $W$ is given by 
$w* \ov{c} = \zeta w$ where $\zeta = e^{2\pi i/N}$. Therefore left action
of $\ov{c} \in \pi$ on $W$ is given by the formula $\ov{c} w = w* \ov{c}$,
and so once again $\ov{c} w = \zeta w$. 
If $\alpha \in Aut(Y/X) - \pi_{x'}$, then $\alpha(W) \cap W = \emptyset$.

\stm {\bf The Grothendieck construction of the bundle $\EE(\sigma)$.}
With notation as in the above lemma,
let $\sigma : \pi \to GL(n)$ be any representation. Then the action of $\pi$ on $Y$ lifts to
an action on the trivial bundle $Y\times \C^n$ over $Y$, where the action on the fibers is
via $\sigma$. Let $\OO^n_Y(\sigma)$ denote this bundle with the given $\pi$-action via $\sigma$.
Let the sheaf 
$$\EE(\sigma) = p^{\pi}_*(\OO^n_Y(\sigma))$$
denote its {\bf invariant direct image on $X$}. The sections of $\EE(\sigma)$ over an open subset
$V\subset X$ are all the sections of $\OO^n_Y$ over $p^{-1}(V)$ that are invariant under $\pi$.
Then $\EE(\sigma)$ is a coherent, torsion free $\OO_X$-module, so it is a vector bundle on $X$.


\stm {\bf The natural holomorphic connection $\nabla$ on $\EE(\sigma)|_{X-x}$.}\\ 
Let $\rho : \pi_1(X-x,y)\to GL(n)$ be a representation with $\rho(c^N) = I$,
so that $\rho$ factors via $\ov{\rho} : \pi_1(X-x,y)/K \to GL(n)$. Let
$\sigma = \ov{\rho}^{op} : \pi  \to GL(n)$.
By statements \ref{ass bun as inv dir im} and \ref{left right correspondence}, 
over $X-x$ the bundle $\EE(\sigma)$
restricts to the vector bundle associated to the principal $\pi_1(X-x,y)/K$-bundle 
$Y'\to X'$ (with base point $y'$), which is the vector bundle 
$E_{\rho}$ on $X-x$ with a frame $y^*$ over $y$ and a connection $\nabla_{\rho}$
as in Section 2. 
If $V\subset X'$ is an open disk, then any flat section
of $E_{\rho}$ over $V$
(means a section of the sheaf $\ker(\nabla)$)
has the form $(\phi, \xi)$ where $\phi : V\to Y'$ is a section, and
$\xi \in \C^n$.

\stm\label{Grothendieck is Deligne}
{\bf Proposition.} {\it Let $\sigma: \pi \to GL(n)$ be such that
  $\sigma(\ov{c}) = e^{2\pi i s/N}I$ where $0\le s < N$, in particular,
  $\sigma(\ov{c}^N) = 0$. If such a $\sigma$ exists then 
  the integer $s$ must be such that $N$ divides $sn$. Let $m = sn/N$, so that 
  $0\le m < n$ and $m/n = s/N$. Let $\rho : \pi_1(X-x,y)\to GL(n)$ be the unique
  representation such that $\sigma = \ov{\rho}^{op} : \pi \to GL(n)$, in particular,
  $\rho(c) = \sigma(\ov{c}^{-1}) = e^{- 2\pi i m/n}I$. 
  Then the Grothendieck bundle $\EE(\sigma)$ is isomorphic as a vector bundle on $X$ with the
  Deligne extension $E(\rho)$ of $E_{\rho}$, with residue $(m/n)I$ at $x$.
  It has degree $-sn/N = -m$.
} 

\proof It is enough to show that natural holomorphic connection $\nabla$ on
$\EE(\sigma)|_{X-x} = E_{\rho}$ extends to a logarithmic connection on
$\EE(\sigma)$ with residue $(s/N)I$ at $x$.
Hence $\deg(\EE(\sigma)) = - sn/N$, and as this has to be an integer, it is necessary 
that $N$ divides $sn$ if $\sigma$ exists.
To see the local behaviour around $x$ of the connection $\nabla$,
consider the disk $U$ around $x$ with coordinate $z$, and the quotient
map $W\to U : w\mapsto z= w^N$ under the action of $\pi_{\wh{x}} = \{ \ov{c},\ldots, \ov{c}^N\}$
with $\ov{c}w = e^{2\pi i /N}w$. 
If $e_i$ is the standard basis of $\C^n$, then 
each $h_i(w) = (w, w^se_i)$ is a $\pi_{\wh{x}}$-invariant holomorphic section of $\OO^n_W(\sigma)$.  
Hence $(h_i)$ is a free basis of holomorphic sections over $U$ for
$\EE(\sigma) = p^{\pi}_*(\OO^n_Y(\sigma))$.
Over a small disk $V\subset U'= U- x$, we have a section $U\to W: z\mapsto z^{1/N}$ of $W\to U$,
for a chosen branch $z^{1/N}$, and so $(z^{1/N}, e_i)$ defines a flat local section.
Note that $h_i = z^{s/N}e_i$. Hence 
$$\nabla(h_i) = \nabla(z^{s/N}e_i) = (s/N)\frac{dz}{z}\otimes z^{s/N}e_i = (s/N)\frac{dz}{z}\otimes h_i$$
hence $\res_x(\nabla) = (s/N)I$ as claimed. 
\hfill$\square$

\stm Conversely, given any $\rho : \pi_1(X-x,y)\to GL(n)$ with $\rho(c) = e^{-2\pi i m/n}$ for some
$0\le m < n$, choose any $N\ge 1$ and $0\le s < N$ such that $s/N = m/n$.
One possible choice is $N=n$ and $s=m$, but there are infinitely many other choices of $N,s$.
Let $p: Y \to X$ be the ramified cover as in Lemma \ref{Poincare ramified cover}.
Let $\sigma = \ov{\rho}^{op} : \pi \to GL(n)$.
Note that $\sigma(\ov{c}) = e^{2\pi i m/n} = e^{2\pi i s/N}$.
Then the above shows that $E(\rho)$
(the Deligne extension with residue $(m/n)I$) is naturally isomorphic to 
the Grothendieck bundle $\EE(\sigma)$.

\section{The cases of genus $0$ and genus $1$}

The existence of Poincar\'e's ramified covering, and hence the
Grothendieck construction which is based on it, needs $g\ge 1$ in general.
That is why it was necessary to assume in [13] that $g\ge 1$.
The assumption in [13] that $g\ge 2$, which is actually unnecessary for the validity of
the main result there (see the argument below), might have come from the Proposition 9.1 in [13],
which indeed needs $g\ge 2$ for it to be true. 

In contrast, the Deligne construction, and the Theorem \ref{main theorem}, work for all $g\ge 0$,
so our reformulation does not put any restriction on the genus. Our proof supplements 
Proposition 9.1 in [13] with the implication (2) $\Rightarrow$ (1) of Lemma
\ref{existence of irr reps genus one} and 
the results of Grothendieck and Atiyah, to cover all $g\ge 0$.

When $g=0$, the group $\pi_1(X-x,y)$ is trivial,
and so there are no irreducible representations of rank $n\ge 2$. 
On the other side, we know from Grothendieck's
decomposition theorem that there are no
stable vector bundles on $X$ of rank $\ge 2$. This reduces the Theorem \ref{main theorem}
to the case $n=1$ when $g=0$. 

When $g=1$, that is, when $X$ is an elliptic curve, Atiyah's theorem [2] implies that 
there exists a stable vector bundle on $X$ of rank $n$ and degree $m$  
if and only if $m$ and $n$ are coprime. To prove the part (2) of Theorem \ref{main theorem}
when $g=1$, we just need to know that there exists an irreducible representation 
$\rho : \pi_1(X-x,y) \to U(n)$ with $\rho(c) = e^{2\pi i m/n}I$ whenever
$m$ and $n$ are coprime. 
Such an example (due to Venkataramana) is given in the proof of
Lemma \ref{existence of irr reps genus one}.

The part (1) of Theorem \ref{main theorem},
together with Atiyah's above result, then proves  
the implication (1) $\Rightarrow$ (2) of the Lemma \ref{existence of irr reps genus one}
(we will not need this implication in the proof of the Theorem \ref{main theorem}). 
The direct elementary proof of (1) $\Rightarrow$ (2) that is given below is due 
to Venkataramana.

\stm\label{existence of irr reps genus one}
{\bf Lemma} {\it Let $n$ and $m$ be integers with $n\ge 1$, and let
$\zeta = e^{2\pi i/n} \in \C$.
Let $F_2$ denote the free group on the two generators $a,b$. 
Let $c = aba^{-1}b^{-1}$.
Then the following two statements are equivalent.\\
$~~$ (1) There exists an irreducible representation $\rho : F_2 \to U(n)$ 
such that $\rho(c) = \zeta^m I$.\\
$~~$ (2) The integers $n$ and $m$ are coprime.

Moreover, regardless of whether $m$ and $n$ are coprime, there
exists an indecomposable representation $\rho : F_2 \to GL(n)$ 
such that $\rho(c) = \zeta^m I$. }

{\bf Proof} (T.N. Venkataramana). \\
(1) $\Rightarrow$ (2) : 
Suppose that $n$ and $m$ are not coprime, in particular, $n\ge 2$.
As $\zeta^n=1$, we can replace $m$ by its residue modulo $n$, so that $0\le m<n$.
We now separately consider the two cases
(a) $m=0$ and (ii) $1 < m < n$.

(i) If $m =0$, then the matrices $\rho(a)=A$ and $\rho(b)=B$ commute, so they are simultaneously
diagonalizable. Hence $\rho$ cannot be irreducible as $n\ge 2$.

(ii) Suppose that $1 < m < n$, and $n$ and $m$ are not coprime.
Let $h \ge 2$ denote the g.c.d. of $n$ and $m$, so that $n = ph$ and $m = qh$,
where $1< p < n$ and $1< q < m$. Note that $mp = nq$.

Suppose that $\rho :F_2\to U(n)$ is an irreducible representation
with $AB = \zeta^m BA$ where $A = \rho(a)$, $B = \rho(b)$. Therefore we have
$$A B^p A^{-1} = (ABA^{-1})^p = \zeta^{mp}B^p = B^p$$
showing that $B^p$ commutes with $A$. As $B^p$ also commutes
with $B$, by irreducibility
%
%
of $\rho$ we must have
$B^p= \lambda I$ for some $\lambda \in \C^{\times}$.

Note that $A$ and $B$ do not commute as by assumption $1< m< n$.
Hence $A$ has an eigenvalue $\alpha$ with multiplicity $<n$.
Let $0\ne v\in \C^n$ such that $Av = \alpha v$.  
Let $V\subset \C^n$ be the subspace
spanned by the vectors $v, Bv, \ldots, B^{p-1}v$. 
As $AB = \zeta^mBA$, for any $k\ge 0$ we have 
$$A B^kv = \zeta^{mk}B^kAv = \zeta^{mk}\alpha B^kv.$$
Hence $V$ is stable under $A$. Also, as $B^p = \lambda I$, $V$ is
stable under $B$. Therefore $V$ is stable under $\rho$,
and $1\le \dim(V) \le p <n$. This contradicts irreducibility of $\rho$.

(2) $\Rightarrow$ (1) : Suppose $n$ and $m$ are coprime. Consider the element $A\in U(n)$
defined in terms of the standard basis $(e_r)$ of $\C^n$ 
by putting $Ae_r = \zeta^{mr} e_r$ for all $1\le r\le n$,
and $Be_1 = e_2,\, \ldots,\, Be_{n-1} = e_n,\, Be_n = e_1$.
Then $A$ has $n$ distinct eigenvalues, so all eigenspaces are $1$-dimensional,
and $B$ cyclically permutes the $n$ eigenspaces of $A$. It follows that $\rho$ is irreducible.
Also, $ABA^{-1}B^{-1} = \zeta^m I$ as required.
(This is an example of a Heisenberg representation of $F_2$ on $U(n)$.)
Hence (1) holds.

Finally, suppose that $n$ and $m$ are not coprime. 
As above, let $h \ge 2$ denote the g.c.d. of $n$ and $m$, so that $n = ph$ and $m = qh$,
where $1< p < n$ and $1< q < m$. Note that $q/p = m/n$, and $p$ and $q$ are coprime.
As above, we have an irreducible unitary representation $\sigma : F_2 \to U(p)$
with $\sigma(a)u_r = e^{2\pi i q/p} u_r$ for all $1\le r\le p$, and
$\sigma(b)u_1 = u_2,\,\ldots,\, \sigma(b)u_p = u_1$ 
where $(u_j)$ is the standard basis of $\C^p$,
so that $\sigma(c) = e^{2\pi i q/p}I = e^{2\pi i m/n}I$ 
(this is a Heisenberg representation of $F_2$ on $U(p)$). 
Now let the (non-unitary) representation $\tau : F_2 \to GL(h)$ be defined 
as follows. We put
$\tau(a) v_1 = v_1,\, \tau(a) v_2 = v_2 + v_1,\,\ldots,\, \tau(a) v_h = v_h + v_{h-1}$
where $(v_j)$ is the standard basis of $\C^h$, and  
we put $\tau(b) = I \in GL(p)$. Then 
for the representation $\rho = \sigma\otimes \tau : F_2 \to GL(\C^p\otimes \C^h) = GL(n)$, 
we have $\rho(c) = e^{2\pi i q/p}I = e^{2\pi i m/n}I$. The only idempotent endomorphisms of $\C^p\otimes \C^h$
that commute with both $\rho(a)$ and $\rho(b)$ are the scalar multiples of $I$, which shows that
$\rho$ is indecomposable. 
\hfill$\square$

\section{A unitary $\rho$ is irreducible $\Leftrightarrow$ $E(\rho)$ is stable}

All the arguments in this section are analogues of the original arguments in [13],
appropriately modified for the Deligne construction in place of the Grothendieck construction,
which makes them simpler.

\subsection{Plurisubharmonic functions}

We recall below some classical facts about harmonic and plurisubharmonic functions on open disks
in $\C$.

\stm\label{plurisubharmonic}
          {\bf Lemma.} {\it Let $U \subset \C$ be an open disk with $0\in U$.
            Then the following statements hold.\\
  (1) If $f$ is a holomorphic function on $U$, then its real and imaginary
  parts are real harmonic functions.\\
  (2) If $\varphi$ is a real harmonic function on the punctured open disk $U - \{ 0 \}$ such that
  $\varphi$ is bounded near $0$, then $\varphi$ uniquely extends to a harmonic function on $U$.\\
%
%
  (3) If $\varphi_1,\ldots,\varphi_n$ are real harmonic functions on $U$, then the
  function $\varphi_1^2 + \cdots + \varphi_n^2$ is plurisubharmonic on $U$.\\
  (4) If $\varphi$ is a real plurisubharmonic function on the open disk $U$ 
  that attains its maximum in $U$, then $f$ is constant.\\
  (5) If $\varphi_1,\ldots,\varphi_n$ are real harmonic functions on $U$ such that 
  $\varphi_1^2 + \cdots + \varphi_n^2$ is constant, then each $\varphi_i$ is constant.
} 
  
\subsection{Proof of Theorem \ref{main theorem}.(1)}

Let $X$ be a compact Riemann surface, and let $x,y\in X$ with $x\ne y$.
Let $n,m$ be any integers, with $n\ge 1$.
Recall from Section 2 that $URep(X-x,y,n,e^{2\pi i m/n})$ denotes the set of all unitary representations
$\rho: \pi_1(X-x,y) \to U(n)$ with local monodromy $\rho(c) = e^{2\pi i m/n}I$ around $x$,
and $ULog(X,x,y,n, - m/n)$ denotes the set of all isomorphism classes of framed
logarithmic connections $(E,\nabla,y^*)$ whose monodromy representation is unitary,
and whose residue at $x$ is $-(m/n)I$. 
By Section 3 we have a bijection $\D_y : URep(X-x,y,n,e^{2\pi i m/n}) \to ULog(X,x,y,n,-m/n)$ defined by
taking the Deligne extension of the connection $(E_{\rho},\nabla_{\rho})$ on $X-x$ for the given
residue $-(m/n)I$, which has the inverse map $\M_y : ULog(X,x,y,n,-m/n) \to 
URep(X-x,y,n,e^{2\pi i m/n})$ given by taking monodromy representations w.r.t. $y^*$.

In particular, when $m=0$, the local monodromy $e^{2\pi i m/n}I$ is $I$, and the residue
$-(m/n)I$ is $0$. Hence for $m=0$ we have the equalities
$URep(X-x,y,n,e^{2\pi i m/n})$ \\
$ = URep(X,y,n)$, the set of all unitary representations
of $\pi_1(X,y)$ of rank $n$, and $ULog(X,x,y,n,-m/n) = UCon(X,y,n)$ the set of all
framed holomorphic connections $(E,\nabla,y^*)$ of rank $n$ on $(X,y)$.

\stm\label{semistability in degree zero}{\bf Lemma on global sections: degree $0$ case.}
{\it With notation as above, let $(E,\nabla, y^*) \in UCon(X,y,n)$, in particular, $\deg(E)=0$.

(a) The inclusion $\ker(\nabla) \hra E$ induces a natural isomorphism
$\Gamma(X,\ker(\nabla)) \to \Gamma(X,E)$.

(b) Composing the above with the natural isomorphism  $(\C^n)^{\rho} \to \Gamma(X,\ker(\nabla))$, 
we have a natural isomorphism $(\C^n)^{\rho} \to \Gamma(X,E)$.

(c) The image of the evaluation map $ev_y: \Gamma(X,E) \to E_y = \C^n$
is the invariant subspace $(\C^n)^{\rho}\subset \C^n$, and the resulting map
$ev_y: \Gamma(X,E) \to (\C^n)^{\rho}$ is an isomorphism,
inverse to the isomorphism in (b) above.
} 

\proof As $E$ is associated to a unitary representation, the standard Hermitian metric on
$\C^n$ (that is preserved by $U(n)$) defines a Hermitian metric on the bundle $E$.
If $s \in \Gamma(W, E)$ a local section over an open subset $W\subset X$, then
the norm square $||s||^2$ is a plurisubharmonic function on $W$. Hence if 
$s \in \Gamma(X, E)$ is a global section then $||s||^2$ is a global plurisubharmonic function
on $X$, so it is constant by Lemma \ref{plurisubharmonic} and the compactness of $X$.
This shows by Lemma \ref{plurisubharmonic} that 
$s$ has constant coefficients w.r.t. any local basis of the unitary local
system $\ker(\nabla)$, so $s$ is a section of $\ker(\nabla)$. Hence the inclusion
$\Gamma(X,\ker(\nabla)) \to \Gamma(X,E)$ is also surjective.
This completes the proof of (a). Now (b) and (c) are clear. \hfill$\square$


\stm\label{the proof of slope inequality for line subundles}
{\bf Lemma on global sections: negative degree case.}
{\it Let the unitary representation
$\rho \in URep(X-x,y,n,e^{2\pi i m/n})$
  be the monodromy representation of a framed logarithmic connection
  $(E,\nabla, y^*)\in ULog(X,x,y,n,-m/n)$, in particular, $\deg(E) = m$.
Then we have the following: If $m < 0$ then $\Gamma(X,E)=0$. }

\proof As $E|_{X-x}$ is associated to the unitary representation $\rho$, which preserves the
standard Hermitian inner product on $\C^n$, it follows that $E|_{X-x}$ has a 
natural hermitian metric over $X-x$. The integrable sections of $E|_{X-x}$
(means local sections of the sheaf $\ker(\nabla|_{X-x})$) 
have constant norms w.r.t. this metric.
By construction of $(E,\nabla,y^*)$ as a Deligne extension of $(E,\nabla,y^*)|_{X-x}$
with residue $-(m/n)I$ at $x$, we have a trivialization
$E|_U \stackrel{\sim}{\to} \OO^n_U$ under which $y^*$ corresponds to the frame $(e_{i,y})$
of $\OO^n_U$ over $y$, where $(e_i)$ is the standard basis of $\OO^n_U$, and
$\nabla|_{U-x}$ has the logarithmic extension
$\nabla_U : \OO^n_U \to \Omega^1_U(\log x) \otimes \OO^n_U$ for which
$$\nabla_U(e_i) = - (m/n) \frac{d z }{z}\otimes e_i.$$
Hence $s_i(z) = z^{m/n}e_i = e^{m/n \log z}e_i$ are (multivalued) integrable sections of 
$\OO^n_{U-x}$, taking the values $e_{i,y}$ over $y$ for the branch of $\log z$ on which $\log 1 =0$.
Hence at any $z\ne 0$, we have 
$|z|^{m/n}||e_i(z)|| = ||s_i(z)|| = ||e_i(y)|| = 1$, which implies that for $z\ne 0$, we have 
$$||e_i(z)|| = |z|^{-m/n}.$$

Now suppose that $\sigma \in \Gamma(X, E)$.
Then the pointwise Hermitian norm $||\sigma ||^2$ is a plurisubharmonic function on $X-x$.
Over $U$, have
$\sigma |_U = \sum\, f_j e_j$ for some holomorphic functions $f_j$ on $U$.
Hence on $U-x$, we have 
$$||\sigma(z)||^2 = (\textstyle\sum_j\, |f_j|^2) |z|^{-2m/n}.$$  

We have assumed that $m<  0$. Hence $\lim_{z\to 0} |z|^{-2m/n}=  0$, and so  
$$\lim_{z\to 0}||\sigma(z)||^2 =  0.$$
Hence by Lemma \ref{plurisubharmonic},
$||\sigma ||^2$ uniquely extends to a plurisubharmonic function on $U$,
which takes the value $0$ at $z=0$.
This shows that $||\sigma||^2$ extends from $X-x$ to $X$
as a plurisubharmonic function $F$ on $X$ with $F(x)=0$. As $X$ is compact, it must attain its maximum.
Hence again by Lemma \ref{plurisubharmonic}, $F$ is identically $0$. Hence $\sigma =0$.
\hfill$\square$

\stm\label{underlying line bundle has smaller slope} {\bf Corollary. } {\it Let
$(E,\nabla, y^*) \in UL(X,x,y,n,-m/n)$ be a logarithmic connection on $(X,x)$
  with a frame $y^*$ for $E_y$ over $y$, such that $\res_x(\nabla) = -(m/n)I$
  (hence with local monodromy $e^{2\pi i m/n }I$ around $x$) and
  such that the monodromy representation $\rho$ for $(E,\nabla, y^*)|_{X-x}$
  is unitary. Let $L\subset E$ be any holomorphic line subbundle. Then
  $\deg(L) \le m/n = \mu(E)$.
} 


\proof Let $p = \deg(L)$, so $\deg(L^*) = -p$ for the dual line bundle $L^*$.
By statement \ref{log con on line bundles} above, $L^*$ has a logarithmic connection $\nabla'$
over $(X,x)$ such that the monodromy of $(L^*, \nabla')$ is unitary, the local monodromy
around $x$ is $1$, and the residue at $x$ is $p$.
Hence the logarithmic connections $\nabla$ and $\nabla'$ induce a logarithmic
connection $\nabla''$ on $E\otimes L^*$ whose monodromy is again unitary,
whose local monodromy around $x$ is $e^{2\pi i m/n }I$ 
and whose residue at $x$ is $(-m/n +p)I$. 
We have a subbundle $\OO_X \subset E\otimes L^*$, so $E\otimes L^*$ has a nowhere
vanishing section. Hence by Lemma \ref{the proof of slope inequality for line subundles},
we must have $m/n - p \ge 0$, which proves the lemma. \hfill$\square$

\stm\label{semistability} {\bf Proposition : Semistability.} 
{\it Let
$(E,\nabla, y^*) \in UL(X,x,y,n,-m/n)$ be a logarithmic connection on $(X,x)$
  with a frame $y^*$ for $E_y$ over $y$, such that $\res_x(\nabla) = -(m/n)I$
  (hence with local monodromy $e^{2\pi i m/n }I$ around $x$) and
  such that the monodromy representation $\rho$ for $(E,\nabla, y^*)|_{X-x}$
  is unitary.
  Let $F\subset E$ be any nonzero holomorphic vector subbundle. Then
  $\deg(F)/\rank(F) \le  \deg(E)/\rank(E)$.
  In other words, the underlying bundle $E$ is a semistable vector bundle on $X$.
} 
          
\proof Let $p = \deg(F)$ and $q = \rank(F)$. Then $L = \wedge^q F$ is a line subbundle
of $E' = \wedge^q E$ of degree $p$. The bundle $\wedge^q E$ has rank
$r = \binom{n}{q}$ and degree $d$ given by
$$d = \deg(\wedge^q E) = \deg(E) \,\binom{\rank(E) -1}{q -1} = m\,\binom{n-1}{q-1}.$$
The logarithmic connection $\nabla$ on $E$ induces
a logarithmic connection $\nabla'$ on $E'$, with 
$$\res_x(E', \nabla') = (-d/r)I$$
(and therefore with local monodromy $e^{2 \pi d/r}I$). 
The Corollary \ref{underlying line bundle has smaller slope} applied to $L \subset E'$
implies that
$p \le d/r$, that is,
$$p \le m\,\binom{n-1}{q-1} / \binom{n}{q} = mq/n.$$
This means $p/q \le m/n$, as desired.
\hfill$\square$

\stm\label{Ramanan}
{\bf The idea of Ramanan.} Let $V$ be a vector space and let $W\subset V$
be a finite dimensional vector subspace.
Consider the $1$-dimensional subspace $\wedge^d W\subset \wedge^d V$ where $d = \dim(W)$.
Let $A : V\to V$ be any linear automorphism such
that under the induced map $\wedge^d A : \wedge^d V \to \wedge^d V $, we have
$(\wedge^d A)(\wedge^d W) = \wedge^d W \subset \wedge^d V$. Then
it can be seen that $A(W) = W \subset V$. Soon after the publication
of [13], Ramanan used the above observation to significantly simplify the proof given in [13] that
the bundles arising from irreducible unitary representations are stable, as it allows
the bypassing of the inductive argument in [13] which uses the full theorem for lower ranks.
We use Ramanan's idea in the proof of (2) $\Rightarrow$ (1) below. 

%
%

\stm\label{stability} {\bf Proposition : Stability.} 
{\it Let
$(E,\nabla, y^*) \in UL(X,x,y,n,-m/n)$ be a logarithmic connection on $(X,x)$
  with a frame $y^*$ for $E_y$ over $y$, such that $\res_x(\nabla) = -(m/n)I$
  (hence with local monodromy $e^{2\pi i m/n }I$ around $x$) and
  such that the monodromy representation $\rho$ for $(E,\nabla, y^*)|_{X-x}$
  is unitary. Then the following two statements are equivalent.\\
(1) The vector bundle $E$ is stable.\\
(2) The representation $\rho$ is irreducible.
} 

\proof (1) $\Rightarrow$ (2) : This argument is essentially the same as
the corresponding argument in [13].
As any unitary representation $\rho$ decomposes into 
direct sum of irreducible unitary representations $\rho_i$'s, the 
framed logarithmic connection $(E,\nabla, y^*)$ is the direct sum of
the corresponding triples $(E_i,\nabla, y^*_i)$.
Note that $\res_x(\nabla)$ is then the direct sum of $\res_x(\nabla_i)$,
so $\res_x(\nabla_i) = (-m/n)I_{n_i}$ where $n_i = \rank(\rho_i) = \rank(E_i)$.
A stable vector bundle is not the direct sum of two nonzero bundles.
Hence if $E_{\rho}$ is stable then $\rho$ must be irreducible. 

(2) $\Rightarrow$ (1) : We will show that if (1) is false then (2) is false. 
We already know from Proposition \ref{stability} that
$E$ is semistable. So if it is not stable, then there would
exists a subbundle $0\ne F\ne  E$ of $E$ with the same slope $p/q = m/n$ as $E$ 
where $p=\deg(F)$, $q=\rank(F)$.
Then $\wedge^q F \subset \wedge^q E$ 
is a line subbundle of degree $p$, hence as proved in Section 4, there exists a unitary representation 
$\theta : \pi_1(X-x, y) \to U(1)$ with trivial local monodromy $1 = e^{2\pi i p}$
(means a unitary representation $\theta : \pi_1(X, y) \to U(1)$) such that
the line bundle $L = \wedge^q F$ admits a logarithmic connection
$\nabla_L$ on $(X,x)$ with monodromy $\theta$ and residue $-p$.
The inclusion $L \hra \wedge^q E$ defines a nowhere vanishing global section
$\sigma \in \Gamma(X, (\wedge^q E)\otimes L^*)$.
The value $\sigma_y$ of $\sigma$ at $y$ corresponds to the inclusion map
$L_y = \wedge^q F_y \hra \wedge^q E_y$.

We have a framed logarithmic triple
$(\wedge^q E\otimes L^*, \nabla', y^*)$ where $\nabla'$ is the
logarithmic connection on $E\otimes L^*$ induced by $\nabla$ and $\nabla_L$.
A simple calculation shows that $\res_x(\nabla')=0$, so $\nabla'$ is actually a holomorphic connection.
By the Lemma \ref{semistability in degree zero}, 
$\sigma \in \Gamma(X,\ker(\nabla')) = \Gamma(X,(\wedge ^qE)\otimes L^*)$. Hence
by Lemma \ref{semistability in degree zero}, $\sigma_y$ is invariant under the monodromy representation
$\eta = (\wedge^q\rho)\otimes \theta^*$ of $\nabla'$, that is,
the inclusion $\wedge^qF_y \hra \wedge^q E_y$ is a $\pi_1(X-x,y)$-equivariant map. 
This means the image of $\sigma_y$, which equals $\wedge^qF_y \subset \wedge^q E$, 
is invariant under $\rho$. 
Hence by \ref{Ramanan} we see that $F_y$ is a $\rho$-invariant subspace of $E_y$.
By assumption $0\ne F_y\ne E_y$, hence $\rho$ is not irreducible, a contradiction. \hfill$\square$

\section{All stable vector bundles arise as $E(\rho)$.}

This section follows the Poincar\'e continuity method of [13], with some simplification
which becomes possible because we use the Deligne construction for the map $\rho \mapsto E(\rho)$ instead
of the Grothendieck construction of the original, and another 
conceptual simplification (well-known since 1970s)
made possible by the Mumford-Seshadri construction of a complex projective moduli variety for 
semistable vector bundles. We sketch these simplifications.
Modulo these, the proof in [13] works verbatim for us, so we do not repeat the
details that are common.

%
%

\newpage

\stm {\bf Table comparing notations.}

\begin{tabular}{ll}

{\bf For us.} & {\bf In [13].}  \\

$UR(X-x, y, n, e^{2\pi i m/n})$ & $U(n,\tau,n)$ \\

$UR^{irr}(X-x, y, n, e^{2\pi i m/n})$ & $U_0(n,\tau,n)$ \\

$UR^{irr}(X-x, n, e^{2\pi i m/n})$ & $M(n,\tau,n)$ \\

$E(\rho)$    &    $p^{\pi}_*(E_{\pi}(\rho))$ \\

$M^s(X,n,m)$ & $S_s(n,m)$ \\

$M^{ss}(X,n,m)$ & (not known in 1965)

\end{tabular}

\stm {\bf Spaces of representations.} 
In particular as shown in [12] and [13], the set of representations 
$UR(X-x, y, n, e^{2\pi i m/n})$ is a closed real analytic 
submanifold of $U(n)^{2g}$, and so it is compact. Moreover, it is 
a Zariski closed real analytic submanifold of the
holomorphic manifold $Rep(X-x,y,n, e^{2\pi i m/n})$,
which is a Zariski closed holomorphic submanifold of $GL(n)^{2g}$. 

\stm\label{domain is nonempty}
{\bf Nonemptyness of the spaces of representations.} The Proposition 9.1 of [13] shows that 
$UR^{irr}(X-x, n, e^{2\pi i m/n})\ne \emptyset$ for all $m,n$ with $n\ge 1$.
We showed in Section 8 that $UR^{irr}(X-x, n, e^{2\pi i m/n})\ne \emptyset$ for
$g=1$ whenever $n$ and $m$ are coprime. Recall that by [2],
stable bundles exist on an elliptic curve only in this case.
For $g=0$, $\pi(X-x) = \{ e \}$, and 
the only case of interest is $n=1$, where the trivial $1$-dimensional representation is irreducible.
Therefore our proof of the Theorem \ref{main theorem} works uniformly for all $g\ge 0$ and $n\ge 1$,
as moreover the Deligne construction $\rho \mapsto E(\rho)$ works for all $g\ge 0$, unlike
the Grothendieck construction in [13] that needs $g\ge 1$ except in special cases.

\stm {\bf Openness of stability.} In [13], the openness
of the stable locus in the parameter space of any family of vector bundles on a curve of $g\ge 2$
was proved in [13] by induction on rank, using the Narasimhan-Seshadri theorem for lower ranks.
A few years later, Narasimhan and Ramanathan [11] 
used the Quot scheme method for proving the openness of the stable  and of
the semistable locus for families of vector bundles on projective varieties of any dimension.
This new method soon became standard, and it also applies widely to families of decorated sheaves
of various kinds.

\stm{\bf Moduli spaces of vector bundles.}
Following Mumford [8] and Seshadri [16], there exists a good moduli variety $M^{ss}(X,n,m)$
of semistable vector bundles on $X$ of rank $n$ and degree $m$, and it is a projective variety.
Its points correspond to S-equivalence classes of semistable bundles, as defined in [16].
This is the same as the set of isomorphism classes of polystable bundles (means direct sums
of stable bundles).
As shown by Narasimhan and Ramanan in [10], the moduli variety $M^s(X,n,m)$ of stable bundles
is a Zariski open dense subvariety of the smooth locus in $M^{ss}(X,n,m)$ (in fact,
Narasimhan and Ramanan [10]
show that $M^s(X,n,m)$ is equal to the smooth
locus of $M^{ss}(X,n,m)$ in all but a few special low rank, low genus cases).
The moduli variety $M^{ss}(X,n,m)$ enjoys the following functorial property.
If $T$ is a real analytic manifold, let $\A_{X\times T}$ be the sheaf of germs of complex valued functions
$f$ on $X\times T$ such that $f$ is real analytic, and for each $t\in T$, the
%
%
restriction $f|_{X\times t}$
is a germ of a holomorphic function on $X\times t$. Let $E_T$ be a locally free 
$\A_{X\times T}$-module of rank $n$ on $X\times T$. Alternatively, we can take
$E_T$ to be a vector bundle of rank $n$ on $X\times T$
defined by an open cover of $X\times T$ with transition functions coming from $\A_{X\times T}$.
Suppose that each restriction $E_t = (E_T)|_{X\times t}$ is semistable of degree $m$.
Then the induced map of sets $T\to M^{ss}(X,n,m)$ that sends $t\mapsto [E_t]$
(the S-equivalence class of $E_t$) is a real analytic map.
If moreover $T$ is a holomorphic manifold and $E$ is a holomorphic bundle, then
the map $T\to M^{ss}(X,n,m): t\mapsto [E_t]$ is holomorphic.
Often in the algebraic geometry literature, this property of the moduli space is stated 
just for algebraic families, but it holds for even continuous families. 
This ultimately follows from the following fact: the universal
property of a complex Grassmannian 
also holds for continuous, $C^{\infty}$, real analytic, holomorphic or algebraic families of 
subbundles of a trivial bundle. From this elementary fact,
it can be concluded by following the construction
of Quot schemes and GIT quotients that the classifying map
$T\to M^{ss}(X,n,m): t\mapsto [E_t]$ of a `nice' family of semistable vector bundles on $X$
is a `nice' map of spaces, where the adjective `nice' can mean any one of
`algebraic', `holomorphic', `real analytic',
`$C^{\infty}$' or just plain `continuous'.

\stm\label{Deligne in families} {\bf Holomorphicity of the Deligne construction in families.}
We have a bijection of sets $\D_y: Con(X-x, y, n, e^{2\pi i m/n}) \to Log(X,x,y,n, -m/n)$ defined
the Deligne construction, with inverse bijection given by the restriction map
${\cal R}_y : Log(X,x,y,n, -m/n) \to Con(X-x, y, n, e^{2\pi i m/n})$, as seen in
Section 3. We claim that by applying the Deligne construction pointwise on the parameter space,
a holomorphic family of relative connections $(E,\nabla)$ on $X$ parameterized by a holomorphic manifold
$T$ gives rise to a holomorphic family of relative logarithmic connections
$\D_y(E,\nabla)$ on $X$ parameterized by $T$ (see [14] for the relevant definitions).
This being a local statement over $T$, we can assume that $T$ is a polydisk in $\C^r$.
We can choose
an open cover $(U_i)$ of $X$ such that each $U_i$ and each $U_{ij} = U_i\cap U_j$ are
biholomorphic to open disks.
As the relative connection $\nabla$ is flat on the
trivial bundle $E|_{U_i\times T}$ over $U_i\times T$,
we can choose corresponding trivializations by flat basic sections.
Let $U_0 = U$ be the open neighbourhood of $x\in X$ chosen at the beginning of
Section 2. Over $U_0\times t$, the logarithmic extension of $\nabla_t$ is described by
$\nabla(e_i) = (-m/n)(dz/z)\otimes e_i$ for the standard basis $e_i$ of $\C^n$. 
This is independent of the point $t\in T$. 
It follows that the family of Deligne extensions will be described by holomorphic
transition functions and holomorphic connection coefficients, so it is a holomorphic
family of logarithmic connections on $X$ parameterized by $T$.

We apply the above to the tautological family of framed holomorphic
connections $({\bf E},\nabla, y^*)$ on $X-x$ parameterized by $T = Con(X-x, y, n, e^{2\pi i m/n})$.
We denote by ${\bf \rho}: \pi_1(X-x,y)\times T \to GL(n)$ the monodromy of this family.
Hence the resulting family ${\bf E}({\bf \rho})$ of vector bundles on $X$
parameterized by $T$, which is the family of 
underlying bundles of the Deligne extension of the family $({\bf E},\nabla, y^*)$, 
is a holomorphic family of vector bundles parameterized by $T$.

\stm\label{proper map}
{\bf The real analytic proper morphisms $\Phi^{ss}: UR(X-x, y, n, e^{2\pi i m/n})\to M^{ss}(X,n,m)$
and $\Phi^s: UR^{irr}(X-x, y, n, e^{2\pi i m/n})\to M^s(X,n,m)$}: 
By Proposition \ref{semistability},
the above holomorphic family of vector bundles on $X$ parameterized by
$Con(X-x, y, n, e^{2\pi i m/n})$,
when restricted to the closed real analytic submanifold
%
%
$UR(X-x, y, n, e^{2\pi i m/n}) \cong UC(X-x,y,n, e^{2\pi i m/n})$, is a real analytic family of
semistable vector bundles on $X$ of rank $n$, degree $m$. Hence its classifying map
is a real analytic morphisms $\Phi^{ss}: UR(X-x, y, n, e^{2\pi i m/n})\to M^{ss}(X,n,m)$ to the
moduli space $M^{ss}(X,n,m)$ of semistable vector bundles.
By Proposition \ref{stability},
the inverse image under $\Phi^{ss}$ of the open subset $M^s(X,n,m)$ (moduli of stable bundles)
is the open subset $UR^{irr}(X-x, y, n, e^{2\pi i m/n})$ of $UR(X-x, y, n, e^{2\pi i m/n})$.
As $UR(X-x, y, n, e^{2\pi i m/n})$ is compact and as $M^{ss}(X,n,m)$ is hausdorff,
the map $\Phi^{ss}$, and therefore its base-change
$\Phi^s: UR^{irr}(X-x, y, n, e^{2\pi i m/n})\to M^s(X,n,m)$, are proper maps of manifolds
(separated and universally closed). 

\stm\label{codomain is connected}
{\bf Connectedness of $M^s(X,n,m)$.} This was proved as Lemma 12.3 in [13].

\stm\label{submersiveness}{\bf $\Phi^s$ is open.} It is proved in [13] that
$\Phi^s$ is a tangent level surjection hence open.
But instead of deducing from this the openness of the induced map
$UR^{irr}(X-x,n,e^{2\pi i m/n}) \to M^s(X,n,m)$ where $UR^{irr}(X-x, n,e^{2\pi i m/n})$ is the quotient
manifold of $UR^{irr}(X-x,y,n,e^{2\pi i m/n})$ by the conjugate action of $U(n)$, they
apply the invariance of domain theorem, leaving out the differential argument. In a conversation
with me fifty years later, Narasimhan ascribed it to their youthful enthusiasm for  
Brouwer's theorem on the invariance of domain.

\stm\label{surjectivity} {\bf Surjectivity of $\Phi^s: UR^{irr}(X-x, y, n, e^{2\pi i m/n})\to M^s(X,n,m)$}: 
The map $\Phi^s$ is continuous, proper, open, its domain is nonempty whenever the codomain
is nonempty, and the codomain is connected. Hence $\Phi^s$ is surjective.
The map $\Phi^s$ factors via $UR^{irr}(X-x, n,e^{2\pi i m/n})$, and the resulting map
is bijective and tangent level isomorphism, so a real analytic isomorphism.
This completes the sketch of the proofs of (2) and (3) in the Main Theorem.

\stm {\bf Generalizations of the logarithmic reformulation.} The Grothendieck construction
of bundles as invariant direct images from representations of the Galois groups of
Poincar\'e's ramified covers have also found use by Ramanathan [15] 
for principal bundles with reductive structure groups
and by Seshadri and by Mehta and Seshadri (see [17] and [7]) for parabolic bundles.
It is possible to apply a similar logarithmic
reformulation to these cases. This will be given elsewhere.

\bigskip

{\bf Acknowledgement. } I thank T.N. Venkataramana for his crucial help involving Heisenberg
representations of the free group $F_2$.

\bigskip

{\footnotesize

{\large \bf References}

[1] Atiyah, M.F. : Complex analytic connections in fiber bundles.
Trans. Amer. Math. Soc. 85 (1957) 181-207.

[2] Atiyah, M.F. : Vector bundles on an elliptic curve.
Proc. London Math. Soc. (3) 7 (1957) 415-452.

[3] Deligne, P. : Equations diff\'erentielles \`a points singuliers r\'eguliers.
LNM 163, Springer Verlag (1970).

[4] Grothendieck, A. : Sur la m\'emoire de Weil ``Generlisations des fonctions ab\'eliennes''.
S\'eminaire Bourbaki, Expose 141 (1956-57).

[5] Kodaira, K. and Spencer, D.C. : On deformations of complex analytic structures, I, II.
Ann. of Math. 67 (1958) 328-466.

[6] Malgrange, B. : Regular connections after Deligne.
Chapter IV in {\it Algebraic $D$-modules} edited by A. Borel, Academic Press, 1987. 

[7] Mehta, V.B. and Seshadri, C.S. : Moduli of vector bundles on curves with parabolic structures.
Math. Annln. 248 (1980) 205-239.

[8] Mumford, D. : {\it Geometric Invariant Theory}. Springer 1965.

[9] Narasimhan, M.S. : Vector bundles on compact Riemann surfaces.
ICTP Trieste lecture notes, IAEA Vienna 1976.

[10] Narasimhan, M.S. and Ramanan, S. : Moduli of vector bundles on a compact Riemann surface.
Ann. of Math. 89 (1969) 14-51.

[11] Narasimhan, M.S. and Ramanathan, A. : Openness of the stability condition on vector bundles.
Unpublished manuscript (1976), TIFR.

[12] Narasimhan, M.S. and Seshadri, C.S. :
Holomorphic vector bundles on a compact Riemann surface.
Math. Annalen 155 (1964) 69-80.

[13] Narasimhan, M.S. and Seshadri, C.S. :
Stable and unitary vector bundles on a compact Riemann surface.
Annals of Math. 82 (1965) 540-567.

[14] Nitsure, N. : Moduli of semistable logarithmic connections. 
Jour. Amer. Math. Soc. 6 (1993) 597-609.

[15] Ramanathan, A. : Stable principal bundles on a compact Riemann surface.
Math. Annln. 213 (1975) 129-152.

[16] Seshadri, C.S. : Space of unitary vector bundles on a compact Riemann surface.
Ann. Math. 85 (1967) 303-336.

[17] Seshadri, C.S. :  Moduli of $\pi$-vector bundles over an algebraic curve.
Questions on Algebraic Varieties, C.I.M.E, III Ciclo, Varenna (1970), p. 139-260.

[18] Weil, A. : Generlisations des fonctions ab\'eliennes, J. Math. Pures et Appl., 17 (1938) 47-87.

\bigskip

Nitin Nitsure,\\
Retired Professor, \\
School of Mathematics,\\
Tata Institute of Fundamental Research,\\
Mumbai 400 005, India. 

email: {\tt nitsure@gmail.com}


}

\end{document}